\documentclass[10pt]{article}

\usepackage[T1]{fontenc}
\usepackage[utf8]{inputenc}
\usepackage{lmodern}
\usepackage[english]{babel}
\usepackage{amsmath}
\usepackage{amssymb}
\usepackage{amsthm}
\usepackage{amsfonts}
\usepackage{stmaryrd}
\usepackage{array}
\usepackage{mathrsfs}
\usepackage{mathtools}
\usepackage[T1]{fontenc}
\usepackage{moresize}
\usepackage{bm}
\usepackage{bbm}
\usepackage{dsfont}
\usepackage[nodisplayskipstretch]{setspace}
\usepackage{fullpage}
\usepackage{lscape}
\usepackage{bigints}
\usepackage{float}
\usepackage{url}
\usepackage{color,xcolor}
\usepackage{multirow}
\usepackage{array}
\usepackage{caption}
\usepackage{hyperref}
\usepackage[pdftex]{graphicx}
\usepackage{booktabs}
\usepackage{longtable}
\usepackage{makecell}
\usepackage[left=25.4mm,right=25.4mm,top=25.4mm,bottom=34.5mm]{geometry}
\usepackage[numbers,sort&compress]{natbib}
\usepackage{sectsty}
\usepackage{enumerate}

\hypersetup{colorlinks=true,allcolors=[rgb]{0,0,1}}

\allowdisplaybreaks

\makeatletter
\DeclareRobustCommand\widecheck[1]{{\mathpalette\@widecheck{#1}}}
\def\@widecheck#1#2{%
    \setbox\z@\hbox{\m@th$#1#2$}%
    \setbox\tw@\hbox{\m@th$#1%
       \widehat{%
          \vrule\@width\z@\@height\ht\z@
          \vrule\@height\z@\@width\wd\z@}$}%
    \dp\tw@-\ht\z@
    \@tempdima\ht\z@ \advance\@tempdima2\ht\tw@ \divide\@tempdima\thr@@
    \setbox\tw@\hbox{%
       \raise\@tempdima\hbox{\scalebox{1}[-1]{\lower\@tempdima\box
\tw@}}}%
    {\ooalign{\box\tw@ \cr \box\z@}}}
\makeatother

\newtheorem{theorem}{Theorem}[section]
\newtheorem{proposition}[theorem]{Proposition}
\newtheorem{lemma}[theorem]{Lemma}
\newtheorem{corollary}[theorem]{Corollary}
\newtheorem{remark}[theorem]{Remark}
\newtheorem{example}[theorem]{Example}

\newcommand{\B}{\mathrm{B}}
\newcommand{\N}{\mathbb{N}}
\newcommand{\R}{\mathbb{R}}
\newcommand{\Z}{\mathbb{Z}}
\newcommand{\E}{\mathbb{E}}
\newcommand{\PP}{\mathbb{P}}

\newcommand{\ind}{\mathds{1}}
\newcommand{\dd}{\,\mathrm{d}}
\newcommand{\Var}{\operatorname{Var}}

\newcommand{\argmax}{\operatorname*{arg\,max}}

\newcommand{\dto}{\xrightarrow{d}}

\newcommand{\trans}{\mathsf{T}}

\DeclareMathOperator{\sgn}{sgn}
\newcommand{\TP}{\mathrm{TP}}
\newcommand{\eff}{\mathrm{eff}}

\begin{document}

\setlength{\abovedisplayskip}{5pt}
\setlength{\belowdisplayskip}{5pt}
\setlength{\abovedisplayshortskip}{5pt}
\setlength{\belowdisplayshortskip}{5pt}
\hypersetup{
  colorlinks,
  breaklinks,
  urlcolor=blue,
  linkcolor=blue,
  pdftitle={Chatterjee's rank correlation under Student-\(t\) dependencee},
  pdfauthor={Sebastian Fuchs}
}

\title{\LARGE Chatterjee's rank correlation under Student-\(t\) dependence}

\author{Sebastian Fuchs\footnote{Department of Artificial Intelligence and Human Interfaces, Paris Lodron Universität Salzburg, Austria, \\ Email: \tt{sebastian.fuchs@plus.ac.at}}}

\date{\today}

\maketitle

\vspace{-1em}
\begin{abstract}
\citet{chatterjee2020} recently introduced a nearest neighbor (NN) estimator of the dependence measure \(\xi\) proposed by \citet{Dette2013}, which quantifies the strength of functional dependence.
Specifically, \(\xi\) equals \(0\) if and only if the variables are independent, and \(1\) if and only if they exhibit perfect functional dependence.
Beyond these characterizations, the behavior of \(\xi\) under specific dependence models remains less well understood.
Moreover, NN-based estimators of \(\xi\) can exhibit substantial bias and high sampling variability, motivating alternative estimation strategies that exploit the underlying dependence structure.
In this paper, we develop a theoretical and inferential framework for \(\xi\) under Student-\(t\) copula dependence.
We derive closed-form expressions for \(\xi\), prove that it is strictly increasing in the absolute value of the correlation parameter and strictly decreasing in the degrees of freedom, and propose a semiparametric estimator that leaves the marginal distributions unspecified.  We establish its strong consistency and asymptotic normality. 
Under correct model specification, our simulations show substantial improvements in estimation accuracy and stability relative to existing estimators. 
By exploiting the Student-\(t\) copula structure, the proposed estimator can also attain an asymptotic variance below the semiparametric efficiency bound recently studied by \citet{DetteHan2026} for the larger model with unrestricted dependence.
\end{abstract}

\section{Introduction} \label{sec:Intro}
Student-$t$ copulas \(C_{\nu,\rho}\), with degrees of freedom \(\nu>0\) and correlation parameter \(|\rho|<1\), are widely used to model dependence, particularly in applications involving joint extreme events.  Examples include flood and drought modeling \citep{nguyen2020biases,cammalleri2024exploring}, wind-power forecasting \citep{kolkmann2024modeling}, portfolio and credit risk assessment \citep{kole2007selecting,frey2003dependent}, and the modeling of dependent insurance claim severities \citep{frees2008hierarchical}. Student-$t$ copulas also serve as building blocks for more complex dependence models, including vine copulas \citep{aas2009paircopula,fuchs2021EoS}, grouped dependence structures \citep{daul2003grouped}, and models with time-varying dependence \citep{kim2016copuladccgarch,cortese2024regimeswitching}.

A distinctive feature of Student-$t$ copulas is that their dependence structure is governed by both \(\rho\) and \(\nu\).  The roles of these parameters can be understood through the Gaussian scale-mixture representation of the underlying bivariate Student-$t$ distribution,
\begin{align} \label{t.ScaleMix}
  (X,Y)
	\overset{d}{=} \sqrt{\frac{\nu}{S_\nu}} \, (Z_1,Z_2),
\end{align}
where \((Z_1,Z_2)\) is a standard bivariate Gaussian random vector with correlation parameter \(\rho\), and \(S_\nu\sim\chi_\nu^2\) is independent of \((Z_1,Z_2)\). While \(\rho\) controls the linear dependence between the Gaussian components, \(\nu\) determines the distribution of the common random scale factor \(\sqrt{\nu/S_\nu}\). Acting as a latent volatility component, this factor simultaneously scales both variables, inducing nonmonotone dependence even when \(\rho=0\).
This dependence at \(\rho=0\) is not captured by certain classical measures of concordance. In particular, Kendall's tau depends only on \(\rho\) and thus satisfies \(\tau(C_{\nu,0})=0\) for all \(\nu>0\). Similarly, Spearman's rho vanishes at \(\rho=0\), although it generally depends on both model parameters. In contrast, mutual information and the tail dependence coefficients capture the influence of \(\nu\), as illustrated in Table~\ref{tab:t_copula_dependence}.
Moreover, Student-$t$ copulas are not stochastically increasing \citep{Rossell2021}, highlighting the nonmonotone nature of their dependence structure.

These properties motivate the use of dependence measures capable of capturing more general forms of functional dependence when evaluating Student-\(t\) dependence. 
A natural candidate is the dependence measure \(\xi\), originally introduced by \citet{Dette2013} and popularized through the rank-based NN estimation procedure proposed by \citet{chatterjee2020}.
For \((X,Y)\) with nondegenerate \(Y\), it is defined by
\begin{align} \label{Tuniv}
  \xi(Y,X)
	& := \frac{\int_{\R} \Var (P(Y \geq y \, | \, X)) \; \mathrm{d} P^{Y}(y)}
					{\int_{\R} \Var (\mathds{1}_{\{Y \geq y\}}) \; \mathrm{d} P^{Y}(y)}\,.
\end{align}
The coefficient takes values in \([0,1]\) with \(\xi(Y,X) = 0\) if and only if \(X\) and \(Y\) are independent, and \(\xi(Y,X) = 1\) if and only if \(Y=f(X)\) almost surely for some measurable function \(f\).
When the marginal distribution functions are continuous and \(C\) denotes the copula of \((X,Y)\), then \eqref{Tuniv} simplifies to 
\begin{align} \label{eq:xipop}
  \xi(Y,X)
  & = 6 \int_{(0,1)^2} (\partial_1 C(u,v))^2 \dd \lambda(u,v) - 2\,,
\end{align}
where \(\partial_1 C\) denotes the partial derivative of \(C\) in the first coordinate, which exists almost everywhere.
Eq.~\eqref{eq:xipop} indicates that \(\xi(Y,X)\), like the other measures of association listed in Table~\ref{tab:t_copula_dependence}, depends only on the copula.
In particular, for a Student-\(t\) copula, \(\xi\) is fully determined by the parameter vector \((\nu,\rho)\).
We therefore use the notation
$$
  \xi(\nu,\rho)
	= \xi(Y,X)\,.
$$
Since the Student-$t$ copula family excludes the independence copula, \(\xi(\nu,\rho)>0\) for all \(\nu>0\) and \(|\rho|<1\).

\begin{table}[t]
\centering
\renewcommand{\arraystretch}{2.3}
\resizebox{\textwidth}{!}{
\begin{tabular}{|c|cccccc|}
\toprule
& \makecell{\textbf{Explicit} \\ \textbf{expression}}
& \makecell{\textbf{The case} \\ {$\rho=0$}}
& \makecell{\textbf{Monotonicity} \\ \textbf{in} $\nu$}
& \makecell{\textbf{Monotonicity} \\ \textbf{in} $\rho$}
& \makecell{\textbf{Dependence} \\ \textbf{concept}}
& \textbf{Reference}
\\
\midrule
\textbf{Kendall's} $\tau$
& $\displaystyle \tfrac{2}{\pi} \arcsin(\rho)$
& $=0$
& constant
& strictly $\nearrow$ in $\rho$
& concordance
& \cite{DemartaMcNeil2005}
\\
\textbf{Spearman's} $\rho_S$\footnotemark
& $\displaystyle \tfrac{6}{\pi} \E \left[\arcsin(\rho V_\nu)\right]$
& $=0$
& \makecell{constant \\ at $\rho=0$}
& strictly $\nearrow$ in $\rho$
& concordance
& \cite{HintzHofertLemieux2020}
\\
\makecell{\textbf{Tail dependence} \\ $\lambda_L=\lambda_U$\footnotemark}
& $\displaystyle 2 \, T_{\nu+1}\bigl(
-\sqrt{\tfrac{(\nu+1)(1-\rho)}{1+\rho}}
\bigr)$
& $\in (0,1/2)$
& strictly $\searrow$ 
& strictly $\nearrow$ in $\rho$
& \makecell{tail \\ dependence}
& \cite{DemartaMcNeil2005}
\\
\makecell{\textbf{Mutual} \\ \textbf{information} $I$\footnotemark}
& $\displaystyle -\tfrac12\log(1-\rho^2)+K(\nu)$
& $\in(0,\infty)$
& strictly $\searrow$ 
& strictly $\nearrow$ in $|\rho|$
& \makecell{information \\ excess}
& \cite{CalsaveriniVicente2009}
\\
\textbf{Chatterjee's} $\xi$
& Theorem~\ref{thm:general}
& \makecell{$\in (0,1/4)$ \\ Eq.~\eqref{cor:zero-copula-correlation}}
& \makecell{strictly $\searrow$ \\ Theorem~\ref{Thm:Monot.Nu}}
& \makecell{strictly $\nearrow$ in $|\rho|$ \\ Theorem~\ref{Thm:MonCCX}}
& \makecell{functional \\ dependence}
& 
\\
\bottomrule
\end{tabular}%
}
\caption{\small Explicit expressions and properties of selected measures of association for the bivariate Student-$t$ copula $C_{\nu,\rho}$, with degrees of freedom $\nu>0$ and correlation parameter $|\rho|<1$.}
\label{tab:t_copula_dependence}
\end{table}
\addtocounter{footnote}{-2}
\footnotetext{\(V_\nu = W_0 / \sqrt{(W_0+W_1)(W_0+W_2)}\) where $W_0,W_1,W_2$ are independent and identically distributed inverse-gamma random variables with shape and scale parameters $\nu/2$.}
\addtocounter{footnote}{1}
\footnotetext{$T_{\nu+1}$ denotes the distribution function of a univariate Student-$t$ distribution with $\nu+1$ degrees of freedom.}
\addtocounter{footnote}{1}
\footnotetext{The function $K$ is given by
$$K(\nu)
  = 2 \, \log \left(\sqrt{\tfrac{\nu}{2\pi}} \, B\left(\tfrac{\nu}{2},\tfrac12\right)\right) 
	    - \tfrac{\nu+2}{\nu}
	    + (\nu+1) \left[\psi\left(\tfrac{\nu+1}{2}\right) - \psi\left(\tfrac{\nu}{2}\right) \right],$$
where $B$ and $\psi$ denote the beta and digamma functions, respectively.}

\pagebreak
Closed-form expressions for \(\xi\) are rare, with known results for the normal distribution and a few other special cases; see, e.g., \cite{ansari2024DeMo,fuchs2025IJAR,Fuchs2024}.
Although the expression in \eqref{eq:xipop} appears simple, deriving explicit representations of \(\xi\) under Student-$t$ dependence is challenging because the partial derivatives involve Student-\(t\) distribution and quantile functions.
Theorem~\ref{thm:general} provides such a representation in terms of generalized hypergeometric functions, which yields simple closed-form expressions for \(\xi(\nu,\rho)\) when the degrees of freedom are integers. In particular, for the Cauchy copula, we obtain
\(\xi(1,\rho) = 1 - 3\sqrt{1-\rho^2}/\pi\) with \(\xi(1,0)\approx 0.04507 > 0\).
Besides making the parameter dependence explicit, these expressions provide exact benchmarks for numerical computations and simulation studies.

Our second contribution concerns the ordering and limiting behavior of \(\xi(\nu,\rho)\).
We establish that \(\xi(\nu,\rho)\) is continuous in both parameters, strictly increasing in \(|\rho|\), and strictly decreasing in \(\nu\).  These monotonicity properties are consistent with the Gaussian scale-mixture representation in
\eqref{t.ScaleMix}.  Increasing \(|\rho|\) strengthens the dependence between the underlying Gaussian components.  On the other hand, since \(\Var(S_\nu/\nu)=2/\nu\), the random variable \(S_\nu/\nu\) converges to one in \(L^2\) as \(\nu\to\infty\). Consequently, the common random scale factor in \eqref{t.ScaleMix} converges to one in probability, and the additional dependence induced by this factor vanishes in the limit.
We further identify the limiting values of \(\xi(\nu,\rho)\) as \(\nu\to\infty\) and \(\nu\to0\), corresponding to the Gaussian copula and a mixture of the Fr{\'e}chet--Hoeffding bounds, respectively. Combined with strict monotonicity in \(\nu\), these limiting results yield sharp lower and upper bounds for \(\xi(\nu,\rho)\). In particular, \(0<\xi(\nu,0)<1/4\) for every \(\nu>0\).

Finally, we propose a semiparametric estimator for \(\xi(\nu,\rho)\) that leaves the marginal distributions unspecified. We adopt the two-step procedure of \citet{ZeeviMashal2002} and \citet{PengWang2014}, which first estimates \(\rho\) from Kendall's tau and then estimates \(\nu\) by maximizing a one-dimensional rank-based pseudo likelihood. 
Inserting these parameter estimates into the derived expression for \(\xi\) yields a plug-in estimator for which we establish strong consistency and asymptotic normality. The resulting estimates are guaranteed to lie in \([0,1]\).
In simulations under correctly specified Student-\(t\) dependence, the plug-in estimator achieves a lower root mean squared error than Chatterjee's original statistic \cite{chatterjee2020}, the multiple NN coefficient proposed by \citet{LinHan2023}, and the semiparametric efficient estimator of \citet{DetteHan2026} in all configurations considered.
The corresponding confidence intervals also exhibit coverage probabilities closer to the nominal level than those of the nonparametric alternatives in most configurations.
The proposed two-step estimator offers two notable advantages. First, it reduces the numerical optimization to a one-dimensional problem, while simulations indicate performance comparable to that of joint maximization.  Second, for most parameter configurations considered, its asymptotic variance is lower than the semiparametric efficiency bound established for general dependence and continuous margins in \cite{KleinLagnouxRochetNguyen2026,DetteHan2026}, with relative variance reductions of up to \(23\%\). 
These findings illustrate the statistical benefits of exploiting a correctly specified Student-\(t\) copula model while retaining the flexibility of unspecified marginal distributions.

The remainder of the paper is organized as follows. Section~\ref{sec:Main} derives the explicit representations, monotonicity properties, and limiting results for Chatterjee's rank correlation under \(t\) dependence. Section~\ref{sec:estimator} develops the two-step plug-in estimator and its asymptotic theory, while Section~\ref{sec:sim} reports the simulation study. 
Proofs of the results in Sections~\ref{sec:Main} and \ref{sec:estimator} are deferred to Sections~\ref{sec:Proof2} and \ref{sec:Proof3}, respectively.

\section{Chatterjee's rank correlation for the Student-\(t\) copula} \label{sec:Main}
Let $T_{\nu,\rho}$ denote the distribution function of the standardized bivariate $t$ distribution with \(\nu\) degrees of freedom and scale matrix 
\(\Sigma_\rho = \left(\begin{smallmatrix}1&\rho\\ \rho&1\end{smallmatrix}\right)\), \(|\rho|<1\).
Further, let $T_{\nu}$, $f_\nu$, and $q_{\nu}=T_{\nu}^{-1}$ denote, respectively, the distribution function, density, and quantile function of the standardized univariate Student-\(t\) distribution.
The Student-\(t\) copula is then given by
\(C_{\nu,\rho}(u,v) = T_{\nu,\rho} (q_{\nu}(u),q_{\nu}(v))\) for all \((u,v) \in [0,1]^2\).

Throughout, we assume that \((X,Y)\) has continuous marginal distribution functions and Student-\(t\) copula \(C_{\nu,\rho}\). Since Chatterjee's rank correlation \(\xi\) is invariant under marginal distributional transformations \cite{ansari2023MFOCI}, i.e.~\(\xi(Y,X) = \xi(F_Y(Y),F_X(X))\), we may, whenever convenient, fix the marginal distributions and consider the standardized bivariate Student-$t$ model $(X,Y) \sim t_\nu(\mathbf{0},\Sigma_\rho)$.
The corresponding conditional distribution is given in \cite{Ding2016} and is restated here for convenience.

\begin{proposition}[Conditional distribution of standardized $t$-distributions]\label{prop:conditional}~
Suppose $(X,Y) \sim t_\nu(\mathbf{0},\Sigma_\rho)$. 
Then \( (Y \mid X=x) \sim t_{\nu+1} (\rho x,s_{\nu,\rho}(x)^2)\), 
where \(s_{\nu,\rho}(x)^2 = (1-\rho^2)\,(\nu+x^2)/(\nu+1)\).
In particular,
$$
  \PP(Y \leq y \mid X=x) 
	= T_{\nu+1}\left( \frac{y-\rho \, x}{s_{\nu,\rho}(x)} \right)\,.
$$
\end{proposition}

\subsection{Explicit representations}
Theorem \ref{thm:general} expresses Chatterjee's rank correlation for the Student-$t$ copula in terms of generalized hypergeometric functions. 
Its derivation relies on the integral representation \(\xi(\nu,\rho) = 1-6I_{\nu,\rho}\) where
\begin{align}\label{eq:productintegral}
 I_{\nu,\rho}
 & = \int_{\R^2} T_{\nu+1}\left( \frac{y-\rho \, x}{s_{\nu,\rho}(x)} \right) \, 
                \left( 1 - T_{\nu+1}\left( \frac{y-\rho \, x}{s_{\nu,\rho}(x)} \right) \right) \, f_\nu(x)f_\nu(y) \dd\lambda(x,y).
\end{align}
Recall that the generalized hypergeometric function \({}_pF_q\) is defined by
$$
 {}_pF_q (a_1,\ldots,a_p; b_1,\ldots,b_q; z)
 =\sum_{k=0}^\infty
 \frac{(a_1)_k\cdots(a_p)_k}{(b_1)_k\cdots(b_q)_k}\frac{z^k}{k!}\,,
$$
where $(a)_0=1$ and $(a)_k=\prod_{j=0}^{k-1}(a+j)$ denotes the rising factorial for \(k \in \N\).
For $\nu\notin\{1/2,3/2,5/2,\ldots\}$ and $z \in (0,1]$, define 
\begin{align*}
 \mathcal A_\nu(z)
 & = {}_4F_3 (\tfrac12,\tfrac12,\tfrac{\nu+1}{2},\tfrac{1-\nu}{2}; 
              \tfrac32,\nu+\tfrac12,\tfrac12-\nu; 
							z),
 \\
 \mathcal B_\nu(z)
 & = {}_4F_3 (\nu+1,\nu+1,\tfrac\nu2+1,\tfrac{3\nu}{2}+1; 
              \nu+\tfrac32,\nu+2,2\nu+1;
							z),
 \\
 e_\nu
 & = \frac{3\cos(\pi\nu/2)\,
 \Gamma(\nu/2+1)\Gamma(3\nu/2+1)\Gamma(-\nu-1/2)}
 {2^{2\nu+1}(\nu+1)\pi^2\Gamma(\nu+1/2)},
\end{align*}
where $\Gamma$ denotes the gamma function. 
The values $\nu\in\{1/2,3/2,5/2,\ldots\}$ are excluded because the third denominator parameter of \(\mathcal A_\nu\) becomes a nonpositive integer, causing the corresponding value \((b_3)_k\) to vanish for sufficiently large \(k\). 
For the admissible parameter values, the defining series of both \(\mathcal A_\nu(z)\) and \(\mathcal B_\nu(z)\) converge for all \(z\in(0,1]\); for details we refer to \cite{DLMF}.

\begin{theorem}[Hypergeometric representation]\label{thm:general}
Suppose that \((X,Y)\) has the unique copula \(C_{\nu,\rho}\) with $\nu>0$ such that $\nu \notin \{1/2,3/2,5/2,\ldots\}$ and \(|\rho| < 1\).
Then 
\begin{align}\label{thm:general:eq}
  \xi(\nu,\rho)
  & = 1-\frac{3\sqrt{1-\rho^2}}{\pi}\mathcal \, \mathcal A_\nu(1-\rho^2) - (1-\rho^2)^{\nu+1} \, e_\nu \, \mathcal B_\nu(1-\rho^2).
\end{align}
For $\nu \in \{1/2,3/2,5/2,\ldots\}$,  
the right-hand side is defined by continuity, that is, as the limit of the entire expression when \(\nu\) approaches the corresponding half-integer value. This limit exists and coincides with \(\xi(\nu,\rho)\). In addition, \eqref{thm:general:eq} extends continuously to \(\rho\in\{-1,1\}\), where \(\xi(\nu,\rho)=1\).
\end{theorem}

Since \(\rho\) only appears in quadratic form in \eqref{thm:general:eq}, 
\(\xi(\nu,-\rho) = \xi(\nu,\rho)\).

The expression derived in Theorem~\ref{thm:general} admits a first simplification when the degrees of freedom are odd integers.\pagebreak

\begin{corollary}[Closed-form expression for odd integer degrees of freedom]\label{cor:odd}
Suppose that \((X,Y)\) has the unique copula \(C_{\nu,\rho}\) with $\nu=2k-1$, $k \in \N$, and \(|\rho| < 1\).
Then 
$$
  \xi(2k-1,\rho)
  = 1-\frac{3\sqrt{1-\rho^2}}{\pi} \; 
	    \sum_{j=0}^{k-1} \frac{(1/2)_j^2 \, (k)_j \, (1-k)_j}{(3/2)_j \, (2k-1/2)_j \, (3/2-2k)_j\,j!} \, (1-\rho^2)^j.
$$
The sum is a polynomial of degree $k-1$ with positive rational coefficients.
The expression extends continuously to \(\rho\in\{-1,1\}\), where \(\xi(2k-1,\rho)=1\).
\end{corollary}

Corollary~\ref{cor:odd} yields the following simple closed-form expressions for \(\nu \in \{1,3\}\).

\begin{example}[Closed-form expressions for \(\xi(\nu,\rho)\) when \(\nu \in \{1,3\}\)]\label{ex:odd}~~
For every \(|\rho| \leq 1\), Chatterjee's rank correlation for the Cauchy copula \(C_{1,\rho}\) and the Student-\(t\) copula \(C_{3,\rho}\) satisfies
$$
  \xi(1,\rho)
   = 1-\frac3\pi\sqrt{1-\rho^2} \qquad \textrm{ and } \qquad
   \xi(3,\rho)
   = 1-\frac{109-4\rho^2}{35\pi}\sqrt{1-\rho^2}\,.
$$
In particular, 
\(\xi(1,0) = 1-\frac3\pi \approx 0.04507 > 0\) and 
\(\xi(3,0) = 1-\frac{109}{35\pi} \approx 0.00869 > 0\).
\end{example}

For even integer degrees of freedom, the expression in Theorem~\ref{thm:general} also reduces to a finite combination of elementary functions, now involving an inverse trigonometric function.

\begin{corollary}[Closed-form expression for even integer degrees of freedom]\label{cor:even}
Suppose that \((X,Y)\) has the unique copula \(C_{\nu,\rho}\) with
\(\nu=2k\), \(k\in\N\), and \(|\rho|<1\). 
Define \(M_0(\rho) = \pi^2/4-(\arcsin|\rho|)^2\) and, for \(j \in \N\), \(M_j\) recursively by
\begin{align}\label{eq:even-M-recursion}
  M_j(\rho)
  & = \frac{2j-1}{2j(1-\rho^2)}M_{j-1}(\rho) - \frac{|\rho|\arcsin|\rho|}{j\sqrt{(1-\rho^2)}} - \frac{1}{2j^2}. 
\end{align}
Then
\begin{align} \label{eq:xi-even}
  \xi(2k,\rho)
  & = 1-\frac{3}{\pi^2} \left[2M_0(\rho) + \sum_{j=k}^{2k-1}c_{k,j}M_j(\rho) + \sum_{\ell=1}^{k}b_{k,\ell}(1-\rho^2)^\ell\right],
\end{align}
where, for \(k\leq j\leq 2k-1\) and \(1\leq \ell \leq k\),
$$ c_{k,j}
   = (-1)^{k+j+1}\frac{k}{j} \binom{j+2k}{k}\binom{k-1}{j-k}
	 \qquad \qquad \qquad
	 b_{k,\ell}
   = \frac{4^\ell}{2\ell^2\binom{2\ell}{\ell}} \frac{\binom{2k-\ell}{k}}{\binom{2k+\ell-1}{k}}.
$$
The expression for \(\xi\) extends continuously to \(\rho\in\{-1,1\}\), where \(\xi(2k,\rho)=1\).
\end{corollary}

Corollary~\ref{cor:even} yields the following closed-form expression for \(\nu = 2\).

\begin{example}[Closed-form expression for \(\xi(\nu,\rho)\) when \(\nu = 2\)]\label{ex:even}
For every \(|\rho|<1\), Chatterjee's rank correlation for the Student-\(t\) copula \(C_{2,\rho}\) equals
$$
  \xi(2,\rho)
  = 1 - \frac{3}{2\pi^2(1-\rho^2)} \left( (1-\rho^2)(4-\rho^2) + (1-4\rho^2) \left[\frac{\pi^2}{4}-(\arcsin\rho)^2\right] + 6\rho\sqrt{1-\rho^2}\,\arcsin\rho \right).
$$
The expression extends continuously to \(\rho\in\{-1,1\}\), where \(\xi(2,\rho)=1\).
In particular, 
\(\xi(2,0) = \frac58-\frac6{\pi^2} \approx 0.01707 > 0\).
\end{example}

\subsection{Monotonicity and continuity in the parameter vector}
As a measure of functional dependence, \(\xi(\nu,\rho)\) is expected to increase with the absolute correlation parameter \(|\rho|\) for fixed \(\nu\), reflecting stronger dependence between the Gaussian components in the scale-mixture representation \eqref{t.ScaleMix}.
The next result establishes this monotonicity property.

\begin{theorem}[Strict monotonicity in the correlation parameter]\label{Thm:MonCCX}
For fixed \(\nu>0\), the mapping $$|\rho| \to \xi(\nu,\rho)$$ is strictly increasing.
\end{theorem}

Since the additional dependence induced by the common random scale factor in \eqref{t.ScaleMix} weakens as \(\nu\) increases and vanishes as \(\nu\to\infty\), it is also natural to expect \(\xi(\nu,\rho)\) to decrease with the degrees of freedom \(\nu\) for fixed \(\rho\). The following result confirms this intuition.

\begin{theorem}[Strict monotonicity in the degrees of freedom]\label{Thm:Monot.Nu}
For fixed $|\rho|<1$, the mapping $$\nu \mapsto \xi(\nu,\rho)$$ is strictly decreasing. 
\end{theorem}

Neither Theorem~\ref{Thm:MonCCX} nor Theorem~\ref{Thm:Monot.Nu} can be directly deduced from \eqref{eq:xipop}, since the first partial derivatives \(\partial_1 C_{\nu,\rho}\) are not pointwise ordered with respect to either \(\nu\) or \(\rho\).

The next result is given in \cite[Corollary 5.3]{AnsariFuchs2026} and states continuity of \(\xi\) in the parameter vector.

\begin{proposition}[Continuity]\label{prop.xi.cont}
If \(\nu_n \to \nu\) and \(\rho_n \to \rho\), then \(\lim_{n \to \infty} \xi(\nu_n,\rho_n) = \xi(\nu,\rho)\).
\end{proposition}

\subsection{Limiting cases} \label{SubSect:LimitCases}

\begin{figure}
  \centering
  \includegraphics[width=0.67\textwidth]{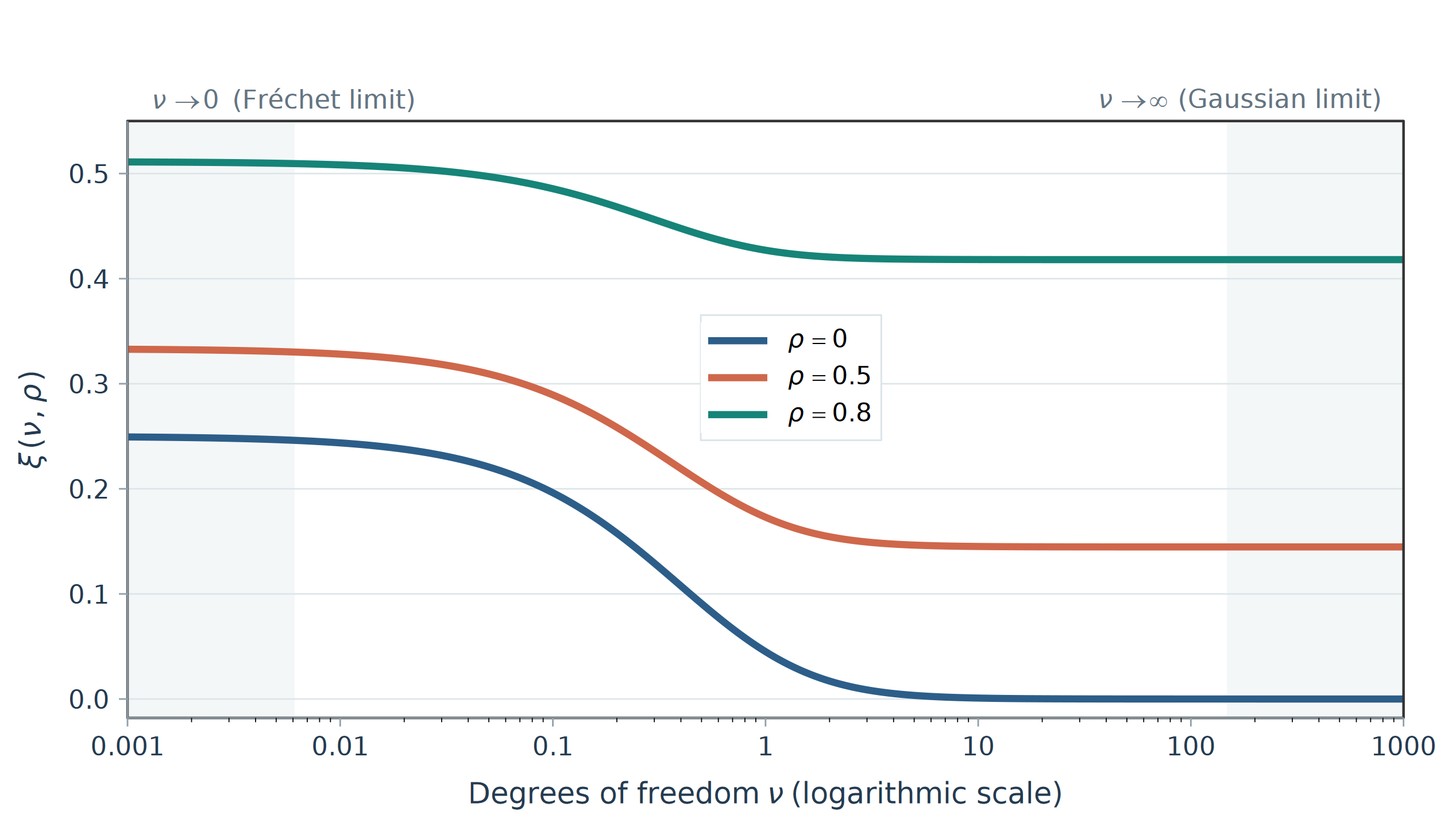}
  \caption{Dependence value \(\xi(\nu,\rho)\) for Student-\(t\) copulas at \(\rho \in \{0,0.5,0.8\}\), as a function of the degrees of freedom \(\nu\).}
  \label{fig:xi-rho-zero}
\end{figure}

While Kendall's tau and Spearman's rho fail to capture the dependence inherent in the \(t\) copula when \(\rho=0\), Chatterjee's rank correlation remains strictly positive and depends on the degrees of freedom.
In this case, Theorem~\ref{thm:general} yields 
\begin{align} \label{cor:zero-copula-correlation}
  \xi(\nu,0) 
	& = 1-\frac{3}{\pi} \; \mathcal A_\nu(1) - e_\nu \, \mathcal B_\nu(1) > 0 
\end{align}
for \(\nu\notin\{1/2,3/2,5/2,\ldots\}\). At the excluded values, the right-hand side is understood as the limit of the entire expression, which exists by continuity. 
Fig.~\ref{fig:xi-rho-zero} illustrates the behavior of \(\xi(\nu,\rho)\) for selected values of \(\rho\) across different values of \(\nu\).

The Student-\(t\) copula \(C_{\nu,\rho}\) converges pointwise to the Gaussian copula \(C_{\rho}\) with correlation parameter \(\rho\) as \(\nu\to\infty\); see \cite{fdsempi2016}. 
For convenience, we write
$$
  \xi(\infty,\rho) 
	:= \frac{3}{\pi}\arcsin\left(\frac{1+\rho^2}{2}\right)-\frac12\,,
$$
for Chatterjee's rank correlation of the Gaussian copula with correlation parameter \(\rho\); see \cite[Example~4]{Fuchs2024}. 
However, because Chatterjee's rank correlation \(\xi\) is not continuous under pointwise convergence of copulas \cite{buecher2024,AnsariFuchs2026}, establishing \(\xi(\nu,\rho) \to \xi(\infty,\rho)\) as \(\nu\to\infty\) requires a stronger notion of convergence.

\begin{theorem}[Gaussian copula limit]\label{Thm:NuToInfty}
For each fixed $|\rho|<1$, as $\nu\to\infty$,
\(\partial_1 C_{\nu,\rho}(u,v) \to \partial_1 C_{\rho}(u,v)\)
almost everywhere on \((0,1)^2\) and in \(L^p\), \(p \in [1,\infty)\). 
As a consequence, 
\(\lim_{\nu \to \infty} \xi(\nu,\rho) = \xi(\infty,\rho)\).
\end{theorem}

Although no proper Student-\(t\) distribution exists at \(\nu=0\), the family admits a meaningful boundary interpretation as \(\nu \to 0\). On the original scale, its mass escapes to \(\pm\infty\). 
As shown in \cite[Theorem 4]{McNeil2022}, the associated sequence of \(t\) copulas converges pointwise to the Fr{\'e}chet mixture
\(C_{0,\rho} := \alpha_\rho \; M + (1-\alpha_\rho) \, W\), where \(\alpha_\rho = \frac12+\frac{\arcsin\rho}{\pi}\).
This limit is a well-defined copula, whose Chatterjee's rank correlation is given by \cite[Example~4]{Fuchs2024}
$$
  \xi(0,\rho)
	:= \frac{1}{4} + \frac{3}{\pi^2} \, \arcsin(|\rho|)^2.
$$
In particular, \(\xi(0,0)= 1/4\).
However, to prove \(\xi(\nu,\rho) \to \xi(0,\rho)\) as \(\nu\to 0\), a stronger notion of convergence is required.

\begin{theorem}[Fr{\'e}chet copula limit]\label{Thm:NuTo0}
For each fixed $|\rho|<1$, as $\nu\to 0$, 
\(\partial_1 C_{\nu,\rho}(u,v) \to \partial_1 C_{0,\rho}(u,v)\)
almost everywhere on \((0,1)^2\) and in \(L^p\), \(p \in [1,\infty)\). 
As a consequence, 
\(\lim_{\nu\to 0}\xi(\nu,\rho) = \xi(0,\rho)\).
\end{theorem}

Theorem~\ref{Thm:NuToInfty} and Theorem~\ref{Thm:NuTo0}, together with the strict monotonicity in \(\nu\) established in Theorem \ref{Thm:Monot.Nu}, yield the following sharp lower and upper bounds for \(\xi(\nu,\rho)\).

\begin{corollary}[Sharp bounds for \(\xi(\nu,\rho)\)]
The strict inequality 
$\xi(\infty,\rho) < \xi(\nu,\rho) < \xi(0,\rho)$
holds for all \(\nu \in (0,\infty)\) and every $|\rho|<1$. 
In particular, \(0 < \xi(\nu,0) < 1/4\).
\end{corollary}

\begin{figure}[t]
  \centering
  \includegraphics[width=0.80\textwidth]{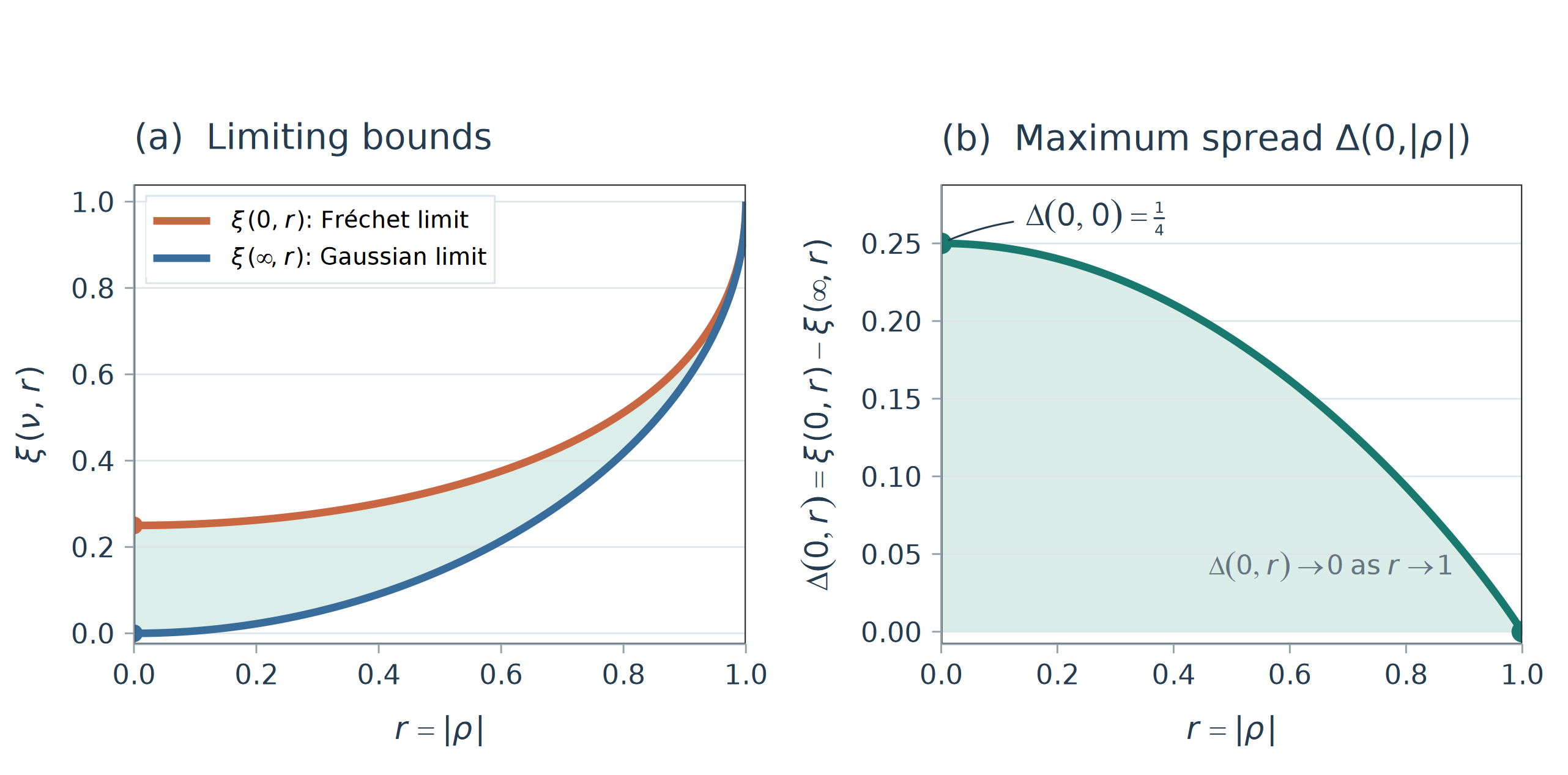}
  \caption{Maximum spread \(\Delta(0,|\rho|)\) for Chatterjee's rank correlation.
(a) Fr{\'e}chet limit $\xi(0,\rho)$ and Gaussian limit $\xi(\infty,\rho)$ as functions of $r=|\rho|$.
(b) The maximum spread $\Delta(0,r)=\xi(0,r)-\xi(\infty,r)$, which is strictly decreasing in $r$.
Its maximum is $\Delta(0,0)=1/4$, while $\Delta(0,r)\to0$ as $r\to 1$.}
  \label{fig:xi-rho-spread}
\end{figure}

Fig.~\ref{fig:xi-rho-spread} depicts the maximum spread \(\Delta(0,|\rho|)\) between the limiting bounds \(\xi(0,\rho)\) and \(\xi(\infty,\rho)\).

\begin{remark}[Excess over the Gaussian copula]
\begin{enumerate}
\item 
For \(\nu \in [0,\infty)\) and \(r=|\rho| \in [0,1]\), define the excess over the Gaussian limit by
$\Delta(\nu,r) := \xi(\nu,r) - \xi(\infty,r)$.
For each fixed \(r < 1\), \(\Delta(\nu,r) > 0\) and, for each fixed \(\nu \in [0,\infty)\), the mapping \(r \to \Delta(\nu,r)\) is strictly decreasing on \([0,1]\) and satisfies \(\Delta(\nu,1) = 0\). 
Moreover,  \(\sup_{\nu \in [0,\infty),\ r\in [0,1]} \Delta(\nu,r) = \frac14 = \xi(0,0)\).

\item 
Mutual information provides a useful contrast. 
For the Student-\(t\) copula with \(\nu \in (0,\infty)\) and \(|\rho| <1\), it admits the additive decomposition \(I(\nu,\rho)=-\frac12\log(1-\rho^2)+K(\nu)\). 
Consequently, its excess over the mutual information of the Gaussian copula with the same correlation parameter is given by \(I(\nu,\rho)-I(\infty,\rho)=K(\nu)>0\), which depends only on the degrees of freedom \(\nu\) and not on \(\rho\).
Thus, changes in \(\rho\) affect mutual information in the same way for all degrees of freedom.
By contrast, the excess \(\Delta(\nu,|\rho|)\) in Chatterjee's rank correlation depends on \(|\rho|\) and converges to zero as \(|\rho|\to1\).
\end{enumerate}
\end{remark}

The strict monotonicity in the degrees of freedom established in Theorem \ref{Thm:Monot.Nu} has the following implication for the identifiability of the parameter vector in the Student-\(t\) copula family.
\pagebreak

\begin{corollary}[Identifiability]\label{cor:df-identifiability}
For each fixed $|\rho|<1$, the mapping \(\nu\mapsto\xi(\nu,\rho)\) is a continuous, strictly decreasing bijection from $(0,\infty)$ onto \(\bigl(\xi(\infty,\rho),\xi(0,\rho)\bigr)\),
with a continuous inverse.
Consequently, the parameter vector $(\nu,\rho)$ of the Student-$t$ copula is uniquely determined by the values of Kendall's tau \(\tau(\nu,\rho)\) and Chatterjee's rank correlation \(\xi(\nu,\rho)\). In particular, the mapping
\((\nu,\rho) \longmapsto \bigl(\tau(\nu,\rho),\xi(\nu,\rho)\bigr)\) is injective on $(0,\infty)\times(-1,1)$.
\end{corollary}

\begin{remark}[Numerical stability]
Numerical evaluation of the \({}_4F_3\) representation in Theorem~\ref{thm:general} may suffer from cancelation errors and require special handling for positive half-integer degrees of freedom.
Likewise, the exact closed-form expression in Corollary~\ref{cor:even} may be prone to cancellation errors, particularly for large even degrees of freedom or correlation parameters near \(\pm1\). 
In Corollary~\ref{prop:joint-one-integral}, we present a one-dimensional integral representation of \(\xi\) that avoids these difficulties and has a nonnegative integrand involving only standard special functions. Valid for all \(\nu>0\), this representation offers a convenient alternative for evaluating the plug-in estimator introduced in Section~\ref{sec:estimator}.
\end{remark}

\section{Estimation}\label{sec:estimator}
We estimate both copula parameters by the two-step procedure proposed by \citet{ZeeviMashal2002} and \citet{PengWang2014}, and then evaluate $\xi$ at the fitted parameters. 
To this end, let \((X_1,Y_1),\ldots,(X_n,Y_n)\) be independent copies of \((X,Y)\), with unknown continuous marginal distribution functions \(F_X,F_Y\) and Student-$t$ copula \(C_{\nu_0,\rho_0}\), where both $\nu_0>0$ and $|\rho_0|<1$ are unknown. 
Recall that the margins of the vector \((X,Y)\) need not be Student-$t$ distributions.
Define $U=F_X(X)$, $V=F_Y(Y)$ and $U_i=F_X(X_i)$, $V_i=F_Y(Y_i)$. 
Then $(U,V)$ and $(U_i,V_i)$ have distribution function $C_{\nu_0,\rho_0}$.

Denote by \(R_i = \sum_{j=1}^n\ind_{\{X_j\leq X_i\}}\) the rank of the predictor and 
by \(S_i = \sum_{j=1}^n\ind_{\{Y_j\leq Y_i\}}\) the response rank.  
Due to continuity of the marginal distribution functions, there are almost surely no ties. 
The pseudo observations and the pseudo log-likelihood are
$$
  \widehat U_i
  = \frac{R_i}{n+1} \qquad\qquad
  \widehat V_i
  = \frac{S_i}{n+1} \qquad\qquad
  L_n(\theta)
  = \frac1n\sum_{i=1}^n \ell_\theta(\widehat U_i,\widehat V_i),
$$
where \(\ell_\theta=\log c_{\nu,\rho}\) and \(c_{\nu,\rho}\) denotes the density of the Student-\(t\) copula \(C_{\nu,\rho}\).  The denominator \(n+1\) ensures that the resulting values lie strictly between \(0\) and \(1\), thereby avoiding the boundary values \(\{0,1\}\) that can lead to difficulties in the evaluation of the quantile functions.

Write \(\theta = (\nu,\rho)\) and \(\theta_0 = (\nu_0,\rho_0)\), and fix a compact rectangle
\(\Theta = [a,b]\times[-r,r]\) such that \((\nu_0,\rho_0) \in (a,b) \times (-r,r)\).
In the first step, define the Kendall's tau statistic and the corresponding estimator of the correlation parameter by
$$
  \widehat\tau_n
	= \binom n2^{-1}\sum_{i<j} \sgn((X_i-X_j)(Y_i-Y_j)) \qquad\qquad
  \widehat\rho_n^{\TP}
  = \sin(\pi\widehat\tau_n/2).
$$
After truncating \(\widehat\rho_n^{\mathrm{TP}}\) to \([-r,r]\), when necessary, choose a measurable maximizer and define
$$
  \widehat\nu_n^{\TP} 
	\in \argmax_{\nu\in[a,b]}L_n(\nu,\widehat\rho_n^{\TP}) \qquad\qquad
	\widehat\theta_n^{\TP} 
	= (\widehat\nu_n^{\TP},\widehat\rho_n^{\TP}) \qquad\qquad
  \xi_n^{\TP}
	= \xi(\widehat\nu_n^{\TP},\widehat\rho_n^{\TP}).
$$
The estimator \(\xi_n^{\TP}\) is semiparametric.  While the copula is characterized by the parameter vector \(\theta\), the marginal distributions remain unspecified and enter the estimation procedure only through the observed ranks. 
Rather than estimating both parameters by joint maximization, the TP procedure estimates \(\rho\) separately and holds this estimate fixed when optimizing over \(\nu\). This reduces the numerical optimization to one dimension, while simulations indicate that the resulting estimators perform comparably to those obtained by joint maximization.

For strong consistency of \(\xi_n^{\TP}\), we require the following uniform convergence and identification result for the pseudo log-likelihood.

\begin{lemma}\label{lem:joint-criterion}
The Student-\(t\) copula family on \(\Theta\) is identifiable. 
The function \(L(\theta)=\E_{\theta_0}[\ell_\theta(U,V)]\) is continuous, has its unique maximum at $\theta_0$, and \(\sup_{\theta\in\Theta}|L_n(\theta)-L(\theta)| \xrightarrow{\mathrm{a.s.}} 0\).
\end{lemma}

The strong consistency of \(\widehat\rho_n^{\TP}\) follows from that of the Kendall's tau estimator, and Lemma~\ref{lem:joint-criterion} yields strong consistency of the maximizer \(\widehat\nu_n^{\TP}\). Together with the continuity of \(\xi\) established in Proposition~\ref{prop.xi.cont}, these results imply strong consistency of the corresponding plug-in estimator of Chatterjee's rank correlation.

\begin{theorem}[Strong consistency]\label{thm:joint-consistency}
We have \(\widehat\theta_n^{\mathrm{TP}} \xrightarrow{\mathrm{a.s.}}\theta_0\) and
\(\xi_n^{\mathrm{TP}}
  \xrightarrow{\mathrm{a.s.}}\xi(\nu_0,\rho_0)\).
\end{theorem}

The joint limit distribution of the two-step estimator \(\widehat\theta_n^{\TP}\) was derived by \citet[Theorem~1]{PengWang2014}. 
We restate their result in Proposition~\ref{prop:tp-peng-wang} and denote the asymptotic covariance matrix by \(\Sigma^{\TP}\); see \cite[Section 5]{PengWang2014} for details.

\begin{proposition}[Limit distribution of \(\widehat\theta_n^{\TP}\)]\label{prop:tp-peng-wang}
We have
\(\sqrt n \, (\widehat\theta_n^{\TP}-\theta_0) \dto N(\mathbf{0},\Sigma^{\TP})\).
\end{proposition}

We now apply the multivariate delta method to derive the asymptotic normality of the plug-in estimator \(\xi_n^{\TP}\) from the asymptotic normality of \(\widehat\theta_n^{\TP}\) as given in Proposition~\ref{prop:tp-peng-wang}.
To this end, we first show differentiability of \(\xi\) with respect to the parameter vector \((\nu,\rho)\).

\begin{lemma}\label{lem:joint-target-derivative}
The mapping \((\nu,\rho) \mapsto \xi(\nu,\rho)\) is continuously differentiable on $(0,\infty)\times(-1,1)$. 
\end{lemma}

Combining Proposition~\ref{prop:tp-peng-wang} and Lemma~\ref{lem:joint-target-derivative} with \cite[Theorem~3.1]{vdV1998} yields asymptotic normality of \(\xi_n^{\TP}\) as follows.

\begin{theorem}[Asymptotic normality]
\label{thm:tp-xi-limit}
We have \(\sqrt n \, (\xi_n^{\TP}-\xi(\nu_0,\rho_0)) \dto N(0,\nabla_\theta\xi(\theta_0)^\trans\Sigma^{\TP}\nabla_\theta\xi(\theta_0))\).
\end{theorem}

In particular, the two-step plug-in estimator is root-$n$ consistent, i.e., 
\(\xi_n^{\mathrm{TP}}-\xi(\nu_0,\rho_0)=O_P(n^{-1/2})\).

\section{Simulation study}\label{sec:sim}

\begin{figure}[t]
  \centering
  \includegraphics[width=0.85\textwidth]{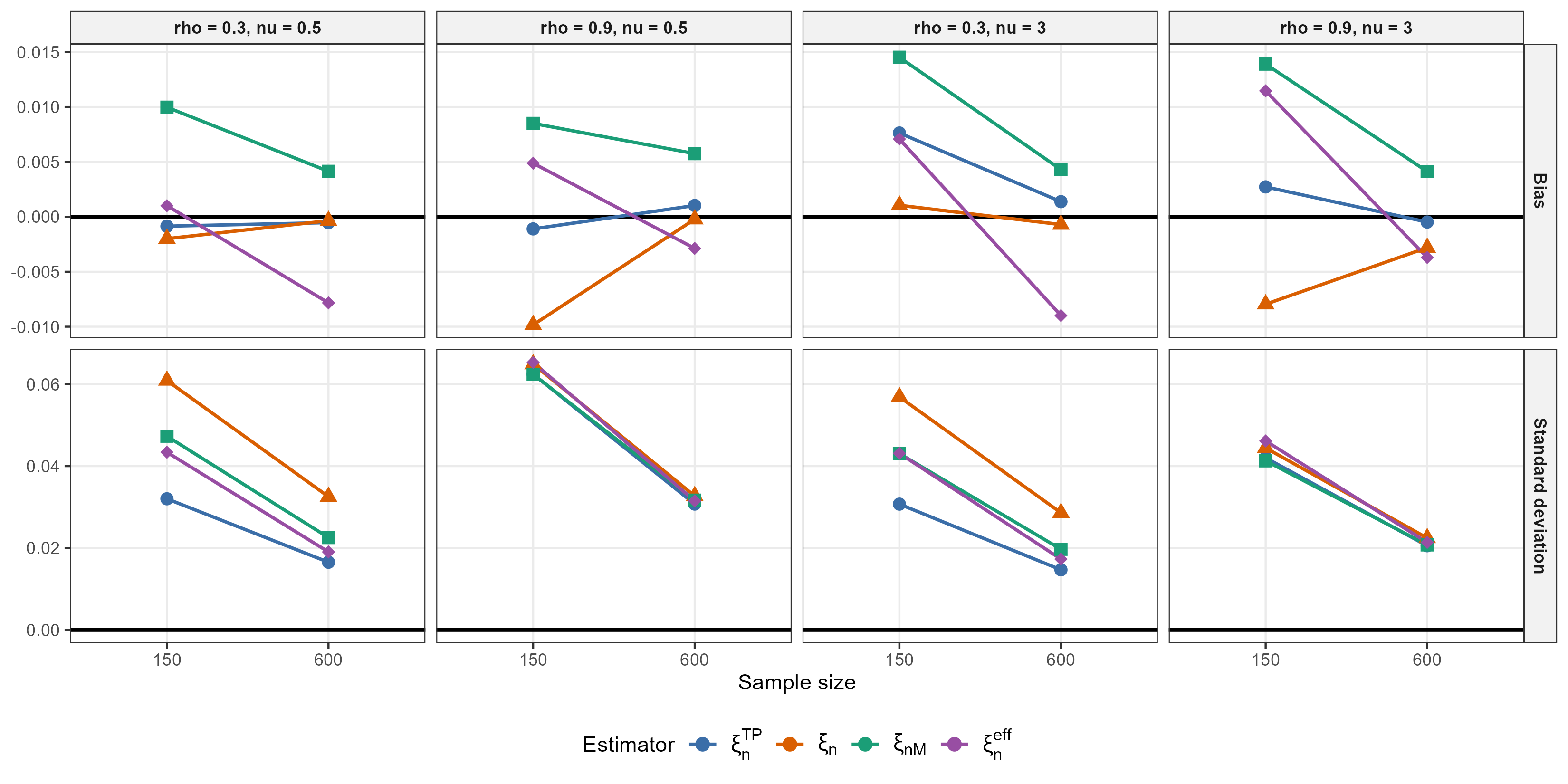}
  \caption{Empirical bias (top row) and standard deviation (bottom row) of the four estimators \(\xi_n^{\mathrm{TP}}\), \(\xi_n\), \(\xi_{n,M}\), and \(\xi_{n}^{\eff}\) under Student-\(t\) dependence for \(n\in\{150,600\}\), \(\rho\in\{0.3,0.9\}\), and \(\nu\in\{0.5,3\}\) over $1,000$ independent repetitions}
  \label{fig:sim.biassd}
\end{figure}

We conduct two simulation studies under Student-\(t\) dependence to compare the two-step plug-in estimator \(\xi_n^{\mathrm{TP}}\) of \(\xi\) introduced in Section~\ref{sec:estimator} with three alternatives:
the original nearest neighbor statistic \(\xi_n\) of \citet{chatterjee2020},
the multiple-right-neighbor estimator \(\xi_{n,M}\) of \citet{LinHan2023}, 
and the semiparametric efficient estimator \(\xi_{n}^{\eff}\) of \citet{DetteHan2026}.
We evaluate the different estimators in terms of bias, standard deviation, root mean squared errors (RMSE), and the empirical coverage probabilities (ECPs) for nominal \(95\%\) confidence intervals.
The bias-corrected nearest neighbor statistic of \citet{han2025Bias} is excluded from the analysis because the required regularity assumptions do not hold simultaneously for the Student-$t$ model, neither on the original scale nor after applying the distributional transform.

We generate independent observations \((X_i, Y_i)\), \(i \in \{1,\dots,n\}\), from the bivariate Student-\(t\)
model
\begin{align*}
  (X,Y)
  & = \sqrt{\frac{\nu}{W}} \, \bigl(Z_{1},\rho Z_{1}+\sqrt{1-\rho^2}\,Z_{2}\bigr),
	& W 
	& \sim \chi^2_\nu
	& Z_{1}, Z_{2}
	& \sim N(0,1),
\end{align*}
where the random variables \(W,Z_1,Z_2\) are independent.
Then \((X,Y)\) is Student-$t$ distributed with $t_\nu$ margins and copula $C_{\nu,\rho}$.
We consider 
\begin{align*}
 n\in\{150,600\},\qquad
 \rho\in\{0,0.3,0.6,0.9\},\qquad
 \nu\in\{0.5,1,3,10\},
\end{align*}
with $R=1000$ independent repetitions for each combination.
The values $\xi(\nu,\rho)$ are computed as described in Section~\ref{sec:Main}, applying Corollary~\ref{cor:odd} for $\nu \in \{1,3\}$, Corollary~\ref{cor:even} for $\nu=10$, and the continuous extension from
Theorem~\ref{thm:general}, evaluated through Corollary~\ref{prop:joint-one-integral}, for $\nu=0.5$.

For \(\xi_n\), we use the R package \texttt{XICOR},  whereas \(\xi_{n,M}\) is computed using the R package \texttt{XIM}; see  \cite{ChenPengZhou2026}.  We set $M_n=\lfloor n^{1/4}\rfloor$, which gives $M_{150}=3$ and $M_{600}=4$. For \(\xi_{n}^{\eff}\), we employ the cross-fitting approach in \cite{DetteHan2026} and estimate the conditional distribution function using \(k_N\)-nearest neighbors on the empirical rank scale, where \(k_N=2\lfloor N^{0.60}/2\rfloor+1\) and \(N=4n/5\) denotes the training sample size.  
For the ECP, the standard deviations of \(\xi_n\) and \(\xi_{n,M}\) are estimated using the \(m\)-out-of-\(n\) bootstrap with \(m=\lfloor n^{1/2} \rfloor\), as recommended by \cite{dette2023boot}.
For \(\xi_n^{\mathrm{eff}}\), we use the empirical standard deviation of the cross-fitted estimates of the efficient influence function.
In contrast, the standard deviations of \(\xi_n^{\mathrm{TP}}\) are obtained from the asymptotic covariance matrix derived in \cite{PengWang2014}, with the required gradient evaluated numerically at the estimated model parameters.

Fig.~\ref{fig:sim.biassd} and Table~\ref{tab.rmse} report the bias, standard deviation, and root mean squared error of \(\xi_n^{\mathrm{TP}}\), \(\xi_n\), \(\xi_{n,M}\), and \(\xi_{n}^{\eff}\) estimated from $1,000$ independent repetitions.
The two-step plug-in estimator \(\xi_n^{\mathrm{TP}}\) exhibits the most favorable bias–standard deviation trade-off and achieves the lowest RMSE across all parameter configurations, outperforming the nonparametric alternatives.
Table~\ref{tab.ecp} reports the empirical coverage probabilities of \(\xi_n^{\mathrm{TP}}\), \(\xi_n\), \(\xi_{n,M}\), and \(\xi_{n}^{\eff}\) estimated from $1,000$ independent repetitions.
The performance of \(\xi_n^{\mathrm{TP}}\) is comparable to that of \(\xi_{n}^{\eff}\) and substantially better than that of \(\xi_n\) and \(\xi_{n,M}\).

\begin{table}[htbp]
\centering
\caption{RMSE for \(\xi_n^{\mathrm{TP}}\), \(\xi_n\), \(\xi_{n,M}\), and \(\xi_{n}^{\eff}\) over $1,000$ independent repetitions}
\label{tab.rmse}
\footnotesize
\setlength{\tabcolsep}{4pt}

\begin{minipage}[h]{0.48\textwidth}
\centering
\begin{tabular}{cccccc}
\toprule
\multicolumn{6}{c}{$\nu = 0.5$} \\
\midrule
$n$ & $\rho$ & $\xi_n^{\mathrm{TP}}$ 
    & $\xi_n$ & $\xi_{n,M}$ & $\xi_{n}^{\eff}$ \\
\midrule
150 & 0   & \textbf{0.0141} & 0.0604 & 0.0394 & 0.0353 \\
600 & 0   & \textbf{0.0061} & 0.0289 & 0.0166 & 0.0149 \\
150 & 0.3 & \textbf{0.0320} & 0.0609 & 0.0483 & 0.0434 \\
600 & 0.3 & \textbf{0.0165} & 0.0325 & 0.0229 & 0.0206 \\
150 & 0.6 & \textbf{0.0577} & 0.0714 & 0.0638 & 0.0636 \\
600 & 0.6 & \textbf{0.0287} & 0.0371 & 0.0314 & 0.0308 \\
150 & 0.9 & \textbf{0.0624} & 0.0656 & 0.0630 & 0.0655 \\
600 & 0.9 & \textbf{0.0307} & 0.0327 & 0.0322 & 0.0316 \\
\midrule
\multicolumn{6}{c}{$\nu = 3$} \\
\midrule
150 & 0   & \textbf{0.0112} & 0.0516 & 0.0328 & 0.0327 \\
600 & 0   & \textbf{0.0037} & 0.0273 & 0.0140 & 0.0130 \\
150 & 0.3 & \textbf{0.0316} & 0.0569 & 0.0454 & 0.0437 \\
600 & 0.3 & \textbf{0.0147} & 0.0286 & 0.0202 & 0.0195 \\
150 & 0.6 & \textbf{0.0492} & 0.0629 & 0.0580 & 0.0572 \\
600 & 0.6 & \textbf{0.0241} & 0.0316 & 0.0268 & 0.0269 \\
150 & 0.9 & \textbf{0.0421} & 0.0451 & 0.0436 & 0.0475 \\
600 & 0.9 & \textbf{0.0206} & 0.0226 & 0.0222 & 0.0217 \\
\bottomrule
\end{tabular}
\end{minipage}
\hfill
\begin{minipage}[t]{0.48\textwidth}
\centering
\begin{tabular}{cccccc}
\toprule
\multicolumn{6}{c}{$\nu = 1$} \\
\midrule
$n$ & $\rho$ & $\xi_n^{\mathrm{TP}}$ 
    & $\xi_n$ & $\xi_{n,M}$ & $\xi_{n}^{\eff}$ \\
\midrule
150 & 0   & \textbf{0.0141} & 0.0547 & 0.0352 & 0.0326 \\
600 & 0   & \textbf{0.0056} & 0.0279 & 0.0154 & 0.0149 \\
150 & 0.3 & \textbf{0.0318} & 0.0584 & 0.0466 & 0.0432 \\
600 & 0.3 & \textbf{0.0160} & 0.0298 & 0.0218 & 0.0213 \\
150 & 0.6 & \textbf{0.0541} & 0.0676 & 0.0626 & 0.0607 \\
600 & 0.6 & \textbf{0.0274} & 0.0344 & 0.0306 & 0.0300 \\
150 & 0.9 & \textbf{0.0548} & 0.0591 & 0.0556 & 0.0581 \\
600 & 0.9 & \textbf{0.0266} & 0.0285 & 0.0273 & 0.0276 \\
\midrule
\multicolumn{6}{c}{$\nu = 10$} \\
\midrule
150 & 0   & \textbf{0.0090} & 0.0531 & 0.0324 & 0.0314 \\
600 & 0   & \textbf{0.0023} & 0.0258 & 0.0139 & 0.0121 \\
150 & 0.3 & \textbf{0.0282} & 0.0553 & 0.0427 & 0.0421 \\
600 & 0.3 & \textbf{0.0135} & 0.0287 & 0.0195 & 0.0183 \\
150 & 0.6 & \textbf{0.0452} & 0.0575 & 0.0533 & 0.0543 \\
600 & 0.6 & \textbf{0.0223} & 0.0295 & 0.0247 & 0.0250 \\
150 & 0.9 & \textbf{0.0379} & 0.0447 & 0.0401 & 0.0446 \\
600 & 0.9 & \textbf{0.0182} & 0.0206 & 0.0191 & 0.0193 \\
\bottomrule
\end{tabular}
\end{minipage}
\end{table}

\begin{table}[htbp]
\centering
\caption{ECPs of nominal $95\%$ intervals for \(\xi_n^{\mathrm{TP}}\), \(\xi_n\), \(\xi_{n,M}\), and \(\xi_{n}^{\eff}\) over $1,000$ independent repetitions.}
\label{tab.ecp}
\footnotesize
\setlength{\tabcolsep}{4pt}

\begin{minipage}[t]{0.48\textwidth}
\centering
\begin{tabular}{cccccc}
\toprule
\multicolumn{6}{c}{$\nu = 0.5$} \\
\midrule
$n$ & $\rho$ & $\xi_n^{\mathrm{TP}}$
    & $\xi_n$ & $\xi_{n,M}$ & $\xi_{n}^{\eff}$ \\
\midrule
600 & 0   & 0.954 & 0.926 & 0.995 & {0.949} \\
600 & 0.3 & 0.914 & 0.920 & 0.980 & {0.935} \\
600 & 0.6 & 0.940 & 0.895 & {0.952} & 0.940 \\
600 & 0.9 & {0.939} & 0.906 & 0.903 & {0.939} \\
\midrule 
\multicolumn{5}{c}{$\nu = 3$} \\
\midrule
600 & 0   & {0.969} & 0.920 & 0.992 & 0.983 \\
600 & 0.3 & {0.942} & 0.926 & 0.978 & 0.931 \\
600 & 0.6 & 0.956 & 0.916 & {0.954} & {0.946} \\
600 & 0.9 & 0.952 & 0.941 & 0.927 & {0.951} \\
\bottomrule
\end{tabular}
\end{minipage}
\hfill
\begin{minipage}[t]{0.48\textwidth}
\centering
\begin{tabular}{cccccc}
\toprule
\multicolumn{6}{c}{$\nu = 1$} \\
\midrule
$n$ & $\rho$ & $\xi_n^{\mathrm{TP}}$
    & $\xi_n$ & $\xi_{n,M}$ & $\xi_{n}^{\eff}$\\
\midrule
600 & 0   & 0.967 & 0.922 & 0.993 & {0.960} \\
600 & 0.3 & {0.945} & 0.926 & 0.984 & 0.930 \\
600 & 0.6 & {0.947} & 0.917 & 0.942 & 0.938 \\
600 & 0.9 & {0.948} & 0.922 & 0.912 & 0.953 \\
\midrule
\multicolumn{5}{c}{$\nu = 10$} \\
\midrule
600 & 0   & 0.990 & {0.937} & 0.992 & 0.985 \\
600 & 0.3 & 0.956 & 0.917 & 0.981 & {0.948} \\
600 & 0.6 & 0.955 & 0.936 & {0.945} & {0.946} \\
600 & 0.9 & 0.943 & 0.946 & 0.925 & {0.949} \\
\bottomrule
\end{tabular}
\end{minipage}
\end{table}

\begin{figure}[htbp]
  \centering
  \includegraphics[width=0.75\textwidth]{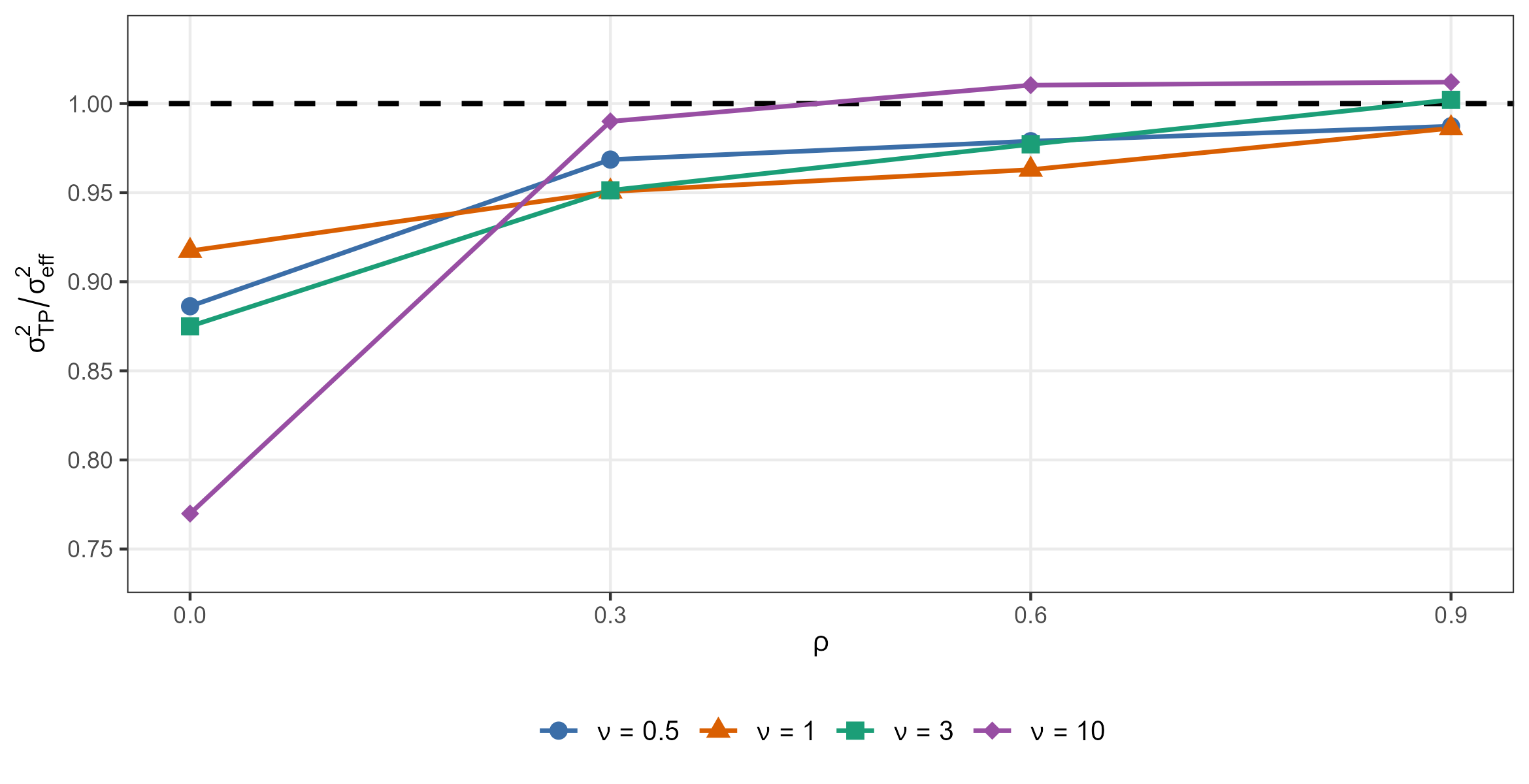}
  \caption{Ratio of the asymptotic variance \(\sigma_{\TP}^2\) of $\xi_n^{\mathrm{TP}}$ to the semiparametric efficiency bound \(\sigma_{\eff}^2\) derived in \cite{DetteHan2026} under Student-$t$ dependence for \(\nu\in\{0.5,1,3,10\}\) and \(\rho\in\{0,0.3,0.6,0.9\}\). Values below the dashed horizontal line indicate a smaller asymptotic variance of the two-step plug-in estimator.}
  \label{fig:sim.lim.var}
\end{figure}

To further investigate the efficiency of the proposed estimator, Fig.~\ref{fig:sim.lim.var} compares its asymptotic variance from Theorem~\ref{thm:tp-xi-limit} with the semiparametric efficiency bound for general dependence and continuous margins derived in \cite{KleinLagnouxRochetNguyen2026,DetteHan2026}. 
The asymptotic variance of \(\xi_n^{\mathrm{TP}}\) is below this bound in $13$ of the $16$ parameter configurations considered, with relative variance reductions of up to $23\%$. In the remaining three configurations, it exceeds the bound only slightly.
These reductions do not contradict the semiparametric efficiency result of \cite{DetteHan2026}, as the latter applies to a substantially broader statistical model with unrestricted dependence.

The proposed two-step plug-in estimator could be further improved by a one-step correction based on the semiparametric efficient score in the Student-$t$ copula model. Simulations suggest that this correction substantially reduces sampling variance, particularly for small degrees of freedom, but may increase finite-sample bias. An additional bias correction could mitigate this effect while retaining the efficiency gains.

\section{Proofs of Section \ref{sec:Main}} \label{sec:Proof2}
In the proofs, we use the following shorthand notation:
\begin{align} \label{Abbrev.Short}
  m 
	& = (\nu+1)/2 \in \big(\tfrac{1}{2},\infty\big)
	& c_m
  & = \frac{1}{\B(\tfrac{1}{2},m)}
	  = \frac{\Gamma(m+1/2)}{\sqrt\pi\,\Gamma(m)} \notag
	\\
	h 
	& = 1-\rho^2 \in (0,1]
	& d_m
	& = \frac{\B(m,m-\frac12)}{\B(m-\frac12,m-\frac12)} 
    = \frac{\Gamma(m)\,\Gamma(2m-1)}{\Gamma(m-1/2)\,\Gamma(2m-1/2)}\,.
\end{align}
\(\Gamma\) denotes the gamma function and we use \(\B(a,b) = \frac{\Gamma(a)\Gamma(b)}{\Gamma(a+b)}\) for the beta function.

We further make use of the following univariate distributions. All densities are zero outside their stated supports.
For the gamma distribution, $b$ is the rate parameter; for the inverse-gamma distribution, $b$ is the scale parameter. 
\(I_x\) denotes the regularized incomplete beta function.

\begin{center}
\renewcommand{\arraystretch}{1.8}
\small
\begin{tabular}{|l|lllll|}
\hline
Distribution & Notation & Parameters & Density $f(x)$ & Support & CDF \\
\hline

Beta
& $\beta(a,b)$
& $a,b>0$
& $\displaystyle \frac{x^{a-1}(1-x)^{b-1}}{\B(a,b)}$
& $x \in (0,1)$
& $I_x$
\\

Beta-prime
& $\beta'(a,b)$
& $a,b>0$
& $\displaystyle \frac{x^{a-1}}{\B(a,b)(1+x)^{a+b}}$
& $x \in (0,\infty)$
& $I_{\tfrac{x}{1+x}}$
\\

Gamma
& $\Gamma(a,b)$
& $a,b>0$
& $\displaystyle \frac{b^a}{\Gamma(a)}x^{a-1}e^{-bx}$
& $x \in (0,\infty)$
&
\\

Inverse-gamma
& $\textrm{Inv}\Gamma(a,b)$
& $a,b>0$
& $\displaystyle \frac{b^a}{\Gamma(a)}x^{-a-1}e^{-b/x}$
& $x \in (0,\infty)$
& 
\\
\hline
\end{tabular}
\end{center}

Before proving Theorem~\ref{thm:general}, we express $I_{\nu,\rho}$, defined in \eqref{eq:productintegral}, in terms of the hypergeometric functions 
\begin{align}
  G_m(\cdot) & := {}_2F_1 \left(\tfrac{1}{2},1-m;m+1;1-\cdot\right), \label{Abbrev.Short.2b}
  \\
  H_m(\cdot) & := {}_2F_1 \left(\tfrac{1}{2},1-m;\tfrac{3}{2};\cdot\right). \label{Abbrev.Short.2}
\end{align}
This allows us to derive the explicit representations of \(\xi\) given in Theorem~\ref{thm:general}, Corollary~\ref{prop:joint-one-integral}, and Lemma~\ref{eq:xi-expectation} from the general representation in \eqref{eq:xipop}.
Throughout, we frequently use hypergeometric identities from \cite{DLMF}.

To derive the desired representations for $I_{\nu,\rho}$, we first rewrite the conditional distribution functions stated in Proposition \ref{prop:conditional} in a more convenient form.

\begin{lemma} \label{App.Lem.Integrand}
For $(X,Y)\sim t_\nu(\mathbf{0},\Sigma_\rho)$, the conditional probability term in \eqref{eq:productintegral} takes the form
$$
  \PP(Y\leq y\mid X=x)\bigl(1-\PP(Y\leq y\mid X=x)\bigr)
  = \frac14 - c_m^2 J_m\left(\frac{y-\rho x}{\sqrt{(y-\rho x)^2+h(\nu+x^2)}}\right)^2,
$$
where \(J_m: [-1,1] \mapsto \R\) is defined by \(J_m(z) = z \, H_m(z^2)\) and with \(m, h, c_m\) as given in \eqref{Abbrev.Short}.
\end{lemma}

\begin{proof} 
Recall that \(I_x\) denotes the regularized incomplete beta function and, for \(|z|\leq1\), define
\begin{align} \label{thm:general:eq4}
  J_m(z)
  & = \sgn(z)\int_{[0,|z|]} (1-t^2)^{m-1}\dd\lambda(t)  
	  = \frac{\operatorname{sgn}(z)}{2} \, B\left( \tfrac{1}{2},m \right) \, I_{z^2}\left( \tfrac{1}{2},m \right)
	  = z \cdot H_m(z^2)\,.
\end{align}
The second identity follows by the substitution \(u=t^2\), while the third follows from \cite[8.17.2 \& 8.17.7]{DLMF}.
The sign factor accounts for negative $z$.
Next, define
\begin{align}\label{Abbrev.Short.3}
  \alpha(x,y)
  & = (y-\rho x)^2+h(\nu+x^2)
    = (x-\rho y)^2+h(\nu+y^2)
		> (y-\rho x)^2, 
  & w
  & = \frac{y-\rho x}{\sqrt{\alpha(x,y)}}
	  \in (-1,1).
\end{align}
The assertion then follows from \(w = \frac{(y-\rho x)/s_{\nu,\rho}(x)}{\sqrt{2m+(y-\rho x)^2/s_{\nu,\rho}(x)^2}}\) and 
\begin{align} \label{thm:general:eq1}
  \PP(Y\leq y\mid X=x)
	& = T_{2m}\left( \frac{y-\rho \, x}{s_{\nu,\rho}(x)} \right)	
		= \frac12+c_mJ_m(w),
\end{align}
where the first identity is due to Proposition~\ref{prop:conditional}, while the second uses symmetry of the $t$ density and \(T_{2m}(z) = \frac12+c_m J_m\left( z/\sqrt{2m+z^2} \right)\). 
\end{proof}

We now prove the announced representations for $I_{\nu,\rho}$.

\begin{lemma} \label{App.Lem.INuRho}
Suppose that \((X,Y)\) has the unique copula \(C_{\nu,\rho}\) with \(\nu>0\) and \(|\rho|<1\). 
Let \(m\), \(h\), and \(d_m\) be as defined in \eqref{Abbrev.Short}, and consider the two independent random variables $B \sim \beta(1/2,m-1/2)$ and $B_m \sim \beta(m,m-1/2)$. 
Then,
\begin{align} 
  I_{\nu,\rho}
	& = \sqrt h \, d_m \, \E[F_m(h B_m)]\,, \label{App.Lem.INuRho.2}
\end{align}
where \(a_m: [0,1) \to \R\) is given by 
\(a_m(z)
  = (1-z)^{-m} \, (1/4 - c_m^2\,z\,H_m(z)^2)\)
and the function \(F_m: (0,1] \to \R\) is defined by
\begin{align}\label{eq:L-integral}
  F_m(s)
  & :=\E[a_m((1-s)B)] 
	   = \frac{c_m^2}{2m\sqrt s} \int_{[0,s]} t^{-1/2} \, G_m(t) \dd\lambda(t),
\end{align}
with \(G_m\) as defined in \eqref{Abbrev.Short.2b}.
\end{lemma}

\begin{proof}
Recall the quantities \(w\) and \(\alpha\) defined in \eqref{Abbrev.Short.3}.
By the invariance property of copulas, we may assume that $(X,Y)\sim t_\nu(\mathbf{0},\Sigma_\rho)$.
Lemma~\ref{App.Lem.Integrand} then yields 
\(\PP(Y\leq y\mid X=x)\bigl(1-\PP(Y\leq y\mid X=x)\bigr) = (1-w^2)^m \, a_m(w^2)\).
Since $1-w^2=h(\nu+x^2)/\alpha(x,y)$ and 
$$
  f_\nu(x)
  = \kappa_m(\nu+x^2)^{-m} \qquad \textrm{with} \qquad
  \kappa_m=\frac{\Gamma(m)\nu^{m-1/2}}
                   {\sqrt\pi\,\Gamma(m-1/2)},
$$
the factors involving $\nu+x^2$ cancel in \eqref{eq:productintegral}. Therefore, 
\begin{align} \label{thm:general:eq2}
  I_{\nu,\rho}
  & = \kappa_m^2 h^m\int_{\R^2}
    \frac{a_m(w^2)}{\alpha(x,y)^m(\nu+y^2)^m}
    \dd\lambda(x,y).
\end{align}
Now, in \eqref{thm:general:eq2}, substitute successively
\(t = (x-\rho y)/\sqrt{h(\nu+y^2)}\), \(a=y/\sqrt{\nu+y^2}\), and \(b=t/\sqrt{1+t^2}\).
The corresponding changes of variables give
\begin{align*}
  \alpha(x,y)
  & = h(\nu+y^2)(1+t^2),
	\\
  (\nu+y^2)^{1/2-2m}\dd\lambda(y)
  & = \nu^{1-2m}(1-a^2)^{2m-2}\dd\lambda(a),
	\\
  (1+t^2)^{-m}\dd\lambda(t)
  & = (1-b^2)^{m-3/2}\dd\lambda(b).
\end{align*}
Using \(w=\sqrt h\,a\sqrt{1-b^2}-\rho b\), \eqref{thm:general:eq2} then reduces to 
$$
  I_{\nu,\rho}
  = \kappa_m^2\sqrt h\,\nu^{1-2m} \int_{(-1,1)^2} a_m\left((\sqrt h\,a\sqrt{1-b^2}-\rho b)^2\right)
      \cdot (1-a^2)^{2m-2}(1-b^2)^{m-3/2} \dd\lambda(a,b).
$$
Since \(\mathrm B(1/2,r) = \int_{(-1,1)}(1-z^2)^{r-1}\dd\lambda(z) = \sqrt{\pi}\,\Gamma(r)/\Gamma(r+1/2)\), normalizing the two weights yields
\begin{align}\label{eq:W}
  I_{\nu,\rho}
  & = \sqrt h\,d_m \; \E[a_m(W^2)],
\end{align}
with \(d_m\) as defined in \eqref{Abbrev.Short}, and where $W=\sqrt h\,A_1\sqrt{1-A_2^2}-\rho A_2$ and $A_1,A_2$ are independent random variables with densities
\begin{align*}
  f_{A_1}(a)
  & = \frac{(1-a^2)^{2m-2}}{\mathrm B(1/2,2m-1)} \; \ind_{[-1,1]}(a),
  & f_{A_2}(b)
  & = \frac{(1-b^2)^{m-3/2}}{\mathrm B(1/2,m-1/2)} \; \ind_{[-1,1]}(b).
\end{align*}

We now replace \(W^2\) in \eqref{eq:W} with an identically distributed random variable. We establish equality in distribution by showing that their moments coincide at every order.
Due to the fact that $A_1$ and $A_2$ are independent and symmetric, only even powers contribute to the binomial expansion
$$
  \E[W^{2k}]
  = \sum_{j=0}^k\binom{2k}{2j} h^j \, (1-h)^{k-j} \, \E[A_1^{2j}] \, \E[A_2^{2(k-j)}(1-A_2^2)^j], \qquad k \in \N.
$$
Since \(A_1^2 \sim \beta(1/2,2m-1)\) and \(A_2^2 \sim \beta(1/2,m-1/2)\), Euler's beta integral gives
\begin{align*}
  \E[A_1^{2j}]
  & = \frac{(1/2)_j}{(2m-1/2)_j},
  & \E[A_2^{2(k-j)}(1-A_2^2)^j]
  & = \frac{(1/2)_{k-j}(m-1/2)_j}{(m)_k}.
\end{align*}
Now, using \(\binom{2k}{2j}(1/2)_j(1/2)_{k-j} = \binom{k}{j}(1/2)_k\) and reindexing the sum yields
$$
  \E[W^{2k}]
    = \frac{(1/2)_k}{(m)_k} \sum_{j=0}^k\binom{k}{j} \, (1-h)^jh^{k-j} \frac{(m-1/2)_{k-j}}{(2m-1/2)_{k-j}}
	  = \frac{(1/2)_k}{(m)_k} \, \E[(1-hB_m)^k]
	  = \E[B^k \, (1-hB_m)^k]
$$
for all \(k \in \N\), with $B\sim\beta(1/2,m-1/2)$ and $B_m\sim\beta(m,m-1/2)$ independent.
Since both $W^2$ and $B(1-hB_m)$ have support $[0,1]$, equality of their moments of all orders implies equality in distribution.
We use this fact in \eqref{eq:reduction} below.
Since $B$ and $B_m$ are independent, \(\E(a_m(B(1-hB_m))\mid B_m) = F_m(hB_m)\) almost surely, and \eqref{eq:W} becomes
\begin{align} \label{eq:reduction}
  I_{\nu,\rho}
	& = \sqrt h \, d_m \,\E[a_m(B(1-hB_m))] 
	  = \sqrt h \, d_m \, \E[F_m(hB_m)]\,,
\end{align}
where we use for the first identity that \(W^2\overset{d}{=}B(1-hB_m)\).
This proves \eqref{App.Lem.INuRho.2}.

Finally, we establish the identity in \eqref{eq:L-integral}.
Let \(B \sim \beta(1/2,m-1/2)\) and define \(f(b)=a_m((1-s)b)\) for \(0<s<1\). Since \(f\) is continuously differentiable on \((0,1)\) and both \(f\) and \(f'\) are bounded on \((0,1]\), the beta Stein identity \cite[Eq.~(2.3)]{GoldsteinReinert2013} yields
$$
  (1-s) \, \E[B(1-B) \, a_m'((1-s)B)]
  = \E\left[\left(mB-\frac12\right) \, a_m((1-s)B)\right]\,.
$$
Differentiating \(a_m\) further gives 
\((1-z) \, a_m'(z)-ma_m(z)
	  = - c_m^2 \, H_m(z)\), using \(H_m(z)+2zH_m'(z)=(1-z)^{m-1}\) \cite[15.5.4 \& 15.4.6]{DLMF}.
Combining these identities with $F_m'(s) = - \E[B \, a_m'((1-s)B)]$ yields
\begin{align} \label{thm:general:eq5}
  & s \, F_m'(s) + \tfrac12 \, F_m(s) 
    = - \E[sB \, a_m'((1-s)B)] + \frac12 \, \E[a_m((1-s)B)] \notag
	\\
	& = - \E[sB \, a_m'((1-s)B)] + m \, \E [B \, a_m((1-s)B)] - (1-s) \, \E[B(1-B) \, a_m'((1-s)B)] \notag
	\\
	& = \E \big[ B \, \big[ m \, a_m((1-s)B)- (1-(1-s)B) \, a_m'((1-s)B) \big] \big] 
	  = c_m^2 \, \E[B \, H_m((1-s)B)].
\end{align}
Expanding the hypergeometric function in \eqref{thm:general:eq5} into its absolutely convergent series and interchanging summation and expectation, together with 
\(\E[B^{k+1}]
  =\frac{(1/2)_{k+1}}{(m)_{k+1}}
  =\frac{1}{2m}\frac{(3/2)_k}{(m+1)_k}\), gives
\begin{align*}
  \E\left[B\,H_m((1-s)B)\right]
	& = \sum_{k=0}^\infty \frac{(1/2)_k(1-m)_k}{(3/2)_k} \frac{(1-s)^k}{k!} \, \E[B^{k+1}]
	\\
  & =
  \frac{1}{2m}\sum_{k=0}^\infty
  \frac{(1/2)_k(1-m)_k}{(m+1)_k}
  \frac{(1-s)^k}{k!}
    = \frac{1}{2m} \, G_m(s).
\end{align*}
Substitution into the preceding identity then yields 
\(s \, F_m'(s)+\tfrac12 \, F_m(s)
  = \frac{c_m^2}{2m}\,G_m(s)\).
Dividing by $\sqrt s$ and applying the product rule then gives
$$
  \frac{\dd}{\dd s}\bigl(\sqrt s\,F_m(s)\bigr)
  = \frac{c_m^2}{2m}s^{-1/2} \, G_m(s).
$$
Since \(\lim_{z\to1} a_m(z) = c_m/(2m)\), \(a_m\) extends continuously to \([0,1]\) and is therefore bounded.
The same holds for $F_m$, implying that $\sqrt s F_m(s)\to0$ as $s\downarrow0$.
Moreover, $G_m(t)$ is bounded near $t=0$ for $m>1/2$. 
Thus the right-hand side is integrable at zero, and integration yields the identity in \eqref{eq:L-integral} which holds for every $m>1/2$.
\end{proof}

\begin{proof}[\textbf{Proof of Theorem~\ref{thm:general}}]
By \eqref{eq:productintegral}, it suffices to evaluate $I_{\nu,\rho}$ using its representation in \eqref{App.Lem.INuRho.2}, \linebreak
$ I_{\nu,\rho}
	= \sqrt h \, d_m \, \E[F_m(h B_m)] $,
with \(F_m\) as given in \eqref{eq:L-integral} and $B_m\sim\beta(m,m-1/2)$.
Substituting \(t=su^2\) in \eqref{eq:L-integral} yields 
\begin{align} \label{eq:L-integral2}
  F_m(s)
  & = \frac{c_m^2}{m} \int_{[0,1]} G_m(su^2) \dd\lambda(u)\,.
\end{align}

Assume first that $\nu\notin\{1/2,3/2,5/2,\ldots\}$ or, equivalently, $2m-1/2\notin \Z$.
For $s\in(0,1)$, \citep[15.10.21]{DLMF} gives
\begin{align*}
  G_m(s) 
	& = {}_2F_1\left(\tfrac{1}{2},1-m;m+1;1-s\right)
	\\
	& = \alpha_m \,{}_2F_1(\tfrac{1}{2},1-m;\tfrac{3}{2}-2m;s)
      + \beta_m s^{2m-1/2} {}_2F_1(m+\tfrac{1}{2},2m;2m+\tfrac{1}{2};s),
\end{align*}
where
\begin{align}\label{Abbrev.Short.4}
  \alpha_m
  & = \frac{\Gamma(m+1)\Gamma(2m-1/2)}
           {\Gamma(m+1/2)\Gamma(2m)}
	\qquad \textrm{and} \qquad
  \beta_m
    = \frac{\Gamma(m+1)\Gamma(1/2-2m)}
           {\sqrt\pi\,\Gamma(1-m)}.
\end{align}
For $m\in\mathbb N$, the coefficient $\beta_m$ vanishes by continuous extension, since $1/\Gamma(1-m)=0$, while the remaining gamma factors are finite.
Applying the substitution \(t=u^2\) and integrating \eqref{eq:L-integral2} term by term, using \citep[Eq.~16.5.2]{DLMF}, yields
\begin{align} \label{eq:L}
  F_m(s)
  & = \frac{c_m^2\alpha_m}{m} \, {}_3F_2 \left(\tfrac{1}{2},\tfrac{1}{2},1\!-\!m;\tfrac{3}{2},\tfrac{3}{2}\!-\!2m;s \right) 
	    + \frac{c_m^2\beta_m s^{2m\!-\!1/2}}{4m^2} \, {}_3F_2 \left(2m,2m,m\!+\!\tfrac{1}{2};2m\!+\!1,2m\!+\!\tfrac{1}{2};s\right).
\end{align}
Now, taking expectations, again using \citep[Eq.~16.5.2]{DLMF}, gives
\begin{align*}
  \E[F_m(hB_m)]
	& = \frac{c_m^2\alpha_m}{m} {}_4F_3 \left(\tfrac{1}{2},\tfrac{1}{2},1\!-\!m,m; \tfrac{3}{2},\tfrac{3}{2}\!-\!2m,2m\!-\!\tfrac{1}{2};h\right)
	\\
  & \quad+\frac{c_m^2\beta_m h^{2m-1/2}}{4m^2} \, \E\left[B_m^{2m-1/2}\right] \, {}_4F_3\left(2m,2m,m+\tfrac{1}{2},3m-\tfrac{1}{2}; 2m+1,2m+\tfrac{1}{2},4m-1;h\right)
	\\
  & = \frac{c_m^2\alpha_m}{m}\,\mathcal A_\nu(h) 
	  + \frac{c_m^2\beta_m}{4m^2} \, h^{2m-1/2} \, \E\left[B_m^{2m-1/2}\right] \, \mathcal B_\nu(h).
\end{align*}
The gamma-function identities
$$
  d_m\frac{c_m^2\alpha_m}{m}
  = \frac1{2\pi}
	\qquad \textrm{ and } \qquad
	6d_m\frac{c_m^2\beta_m}{4m^2} \, \E\left[B_m^{2m-1/2}\right]
  = e_\nu
$$
now yield \eqref{thm:general:eq} for $0<h<1$.

Both ${}_4F_3$ series have parameter excess $1/2$ and therefore converge absolutely at $h=1$ for the degrees of freedom considered here. By Proposition~\ref{prop.xi.cont}, continuity of \eqref{thm:general:eq} in $\rho$ and $\nu$ extends the identity to $h=1$ (that is, $\rho=0$). Moreover, the right-hand side of \eqref{thm:general:eq} extends to positive half-integer degrees of freedom, although its individual terms may be singular at these values.

Finally, as $\rho$ approaches either endpoint \(\{-1,1\}\), the \(t\) copula converges weakly to the corresponding Fr{\'e}chet-Hoeffding bound, \(M\) or \(W\), both of which have Chatterjee's rank correlation equal to \(1\).
By the lower semicontinuity of $\xi$ under weak convergence \cite[Theorem~2.14]{Borgonovo2025}, combined with the upper bound \(\xi(\nu,\rho)\leq1\), the claimed boundary values follow.
\end{proof}

The following explicit expression, based on Lemma~\ref{App.Lem.INuRho} and \eqref{eq:L-integral2}, represents \(\xi(\nu,\rho)\) as an expectation and is used in the proofs of Corollary~\ref{cor:odd}, Theorem~\ref{Thm:MonCCX}, Theorem~\ref{Thm:Monot.Nu}, and Lemma~\ref{lem:joint-target-derivative}.

\begin{lemma}\label{eq:xi-expectation}
Suppose that \((X,Y)\) has the unique copula \(C_{\nu,\rho}\) with \(\nu>0\) and \(|\rho|<1\). 
Let \(m\) and \(h\) be as defined in \eqref{Abbrev.Short}. 
Then, 
$$
  \xi(\nu,\rho)
  = 1-\frac{3\sqrt h}{\pi} \int_{[0,1]} \E[(1+hu^2W_m)^{m-1}] \dd \lambda(u)\,,
$$
where \(W_m=B_m D_m\), with independent random variables \(B_m \sim\beta\left(m,m-\frac12\right)\) and \(D_m \sim \beta'\left(\frac12,2m-\frac12\right)\).
\end{lemma}

\begin{proof}
By \eqref{App.Lem.INuRho.2} and \eqref{eq:L-integral2},
\begin{align} \label{eq:reduction2}
  I_{\nu,\rho}
	&  = \frac{c_m^2\,d_m\,\sqrt{h}}{m} \int_{[0,1]} \E[G_m(h B_m u^2)] \dd\lambda(u)\,.
\end{align}
Euler's integral representation \cite[Eq.~15.6.2\_5]{DLMF} gives
\begin{align*} 
  G_m(s)
	& = \frac{\Gamma(m+1)}{\Gamma(\frac12)\Gamma(m+\frac12)} 
      \, \int_{(0,\infty)} t^{-1/2}(1+t)^{-2m}(1+st)^{m-1} \dd \lambda(t)
	\\
	& = \frac{\Gamma(m+1)}{\Gamma(\frac12)\Gamma(m+\frac12)} 
      \, \B \left(\tfrac12,2m-\tfrac12\right) \, \E\left[(1+sD_m)^{m-1}\right]
		= \alpha_m \, \E\left[(1+sD_m)^{m-1}\right],
\end{align*}
with \(\alpha_m\) as defined in \eqref{Abbrev.Short.4} and where \(D_m \sim \beta' \left(\frac12,2m-\frac12\right)\).
The integral representation is valid for $m>\frac12$ and $s>0$, and extends to $s=0$ by continuity.
Choosing $B_m$ and $D_m$ independently and setting \(W_m=B_mD_m\), \eqref{eq:reduction2} reduces to
$$
	I_{\nu,\rho}
	  = \frac{c_m^2\,d_m\,\alpha_m\,\sqrt{h}}{m} \int_{[0,1]} \E\left[(1+hu^2W_m)^{m-1}\right] \dd \lambda(u)\,.
$$
Since \(c_m^2 d_m\alpha_m/m = 1/2\pi\), the assertion follows from \(\xi=1-6I_{\nu,\rho}\).
\end{proof}

\begin{corollary}[One-dimensional integral representation]\label{prop:joint-one-integral}~
Suppose that \((X,Y)\) has the unique copula \(C_{\nu,\rho}\) with \(\nu>0\) and \(|\rho|<1\). 
Let \(m\) and \(h\) be as defined in \eqref{Abbrev.Short}. 
Then,
$$
  \xi(\nu,\rho)
    = 1-\frac{6c_m^2\sqrt h}{m} \, 
      \int_{[0,1]} G_m(hz^2) \left(1-I_{z^2}\left(m-\tfrac{1}{2},m-\tfrac{1}{2}\right))\right)
      \dd\lambda(z)\,,
$$
where \(I_x\) denotes the regularized incomplete beta function.
\end{corollary}

\begin{proof}
The substitution \(z=u\sqrt{B_m}\) in \eqref{eq:reduction2} yields
$$
  I_{\nu,\rho}
    = \frac{c_m^2\,d_m\,\sqrt h}{m} \, \int_{[0,1]} G_m(hz^2) \; 
      \E\left[ B_m^{-1/2} \, \ind_{\{B_m\geq z^2\}} \right] \dd\lambda(z)\,,
$$
where $B_m \sim \beta(m,m-1/2)$.
For \(b\in(0,1)\), the density of \(B_m\) satisfies
$$
  b^{-1/2}f_{B_m}(b)
  =
  \frac{1}{d_m}
  \frac{b^{m-3/2}(1-b)^{m-3/2}}
       {\mathrm B(m-1/2,m-1/2)}\,.
$$
Consequently,
\(\E\big[B_m^{-1/2}\ind_{\{B_m\geq z^2\}}\big]
  = (1-I_{z^2}(m-\tfrac{1}{2},m-\tfrac{1}{2}))/d_m\).
Substituting this expression into the integral above and recalling that \(\xi=1-6I_{\nu,\rho}\) completes the proof.
\end{proof}

For later use, note that independence of \(B_m \sim\beta\left(m,m-\frac12\right)\) and \(D_m \sim \beta'\left(\frac12,2m-\frac12\right)\) in Lemma~\ref{eq:xi-expectation} yields
\begin{align}\label{eq:W_moments}
  \E[W_m^k]
 	& = \E[B_m^k] \, \E[D_m^k]
    = \frac{(m)_k} {(2m-\frac12)_k} \, \frac{(\frac12)_k}{\prod_{j=0}^{k-1}(2m-\frac32-j)}\notag
  \\
	& = \frac{(-1)^k(m)_k(\frac12)_k}{(2m-\frac12)_k(\frac32-2m)_k} 
	  = \left(\frac12\right)_k \prod_{j=0}^{k-1} \frac{m+j}{4m(m-1)+\frac34-j(j+1)} < \infty
\end{align}
for every $k \in \N_0$ with $k<2m-\frac12$.

\begin{proof}[\textbf{Proof of Corollary~\ref{cor:odd}}] 
Applying the binomial theorem to the expression in Lemma~\ref{eq:xi-expectation} and using \eqref{eq:W_moments} yields 
$$
  \xi(2k-1,\rho)
    = 1 - \frac{3\sqrt h}{\pi} \sum_{j=0}^{k-1} \binom{k-1}{j} \frac{h^j \, \E[W_k^j]}{2j+1}
	  = 1 - \frac{3\sqrt h}{\pi} \sum_{j=0}^{k-1} \binom{k-1}{j} \frac{h^j}{2j+1} \, \frac{(k)_j(\frac12)_j(-1)^j}{(2k-\frac12)_j(\frac32-2k)_j}\,.
$$
The assertion now follows from the identities \(\binom{k-1}{j}=(-1)^j(1-k)_j /j!\) and \(\frac1{2j+1}=(\frac12)_j / (\frac32)_j\).
\end{proof}

\begin{proof}[\textbf{Proof of Corollary~\ref{cor:even}}] 
We use the hypergeometric representation in Theorem~\ref{thm:general}.
Throughout this proof, let \(\nu=2k\), \(k\in\N\), and \(h=1-\rho^2\) as defined in \eqref{Abbrev.Short}. Consider \(0<h<1\).

We first establish the recursion in \eqref{eq:even-M-recursion}.
For \(j\in\N_0\), define
\begin{align}\label{eq:even-readable-integral}
 M_j(\rho)
 & = h^{-j} \int_{[0,h]} \frac{t^{j-1/2}\arccos\sqrt t}{\sqrt{1-t}} \dd\lambda(t).
\end{align}
For \(j=0\), the substitution \(t = (\cos u)^2\) gives
\(M_0(\rho)
  =\left(\frac\pi2\right)^2-(\arccos\sqrt h)^2
   =\frac{\pi^2}{4}-(\arcsin|\rho|)^2\).
For \(j \in \N\), define \(g(t)=\arccos\sqrt t/\sqrt{t(1-t)}\) and \(Q_j(t) = t^{j-1/2} \, \sqrt{1-t} \, \arccos(\sqrt t)\).
By definition,
\begin{align}\label{eq:even-readable-integral2}
  \int_{[0,h]} t^jg(t) \dd\lambda(t) = h^j\,M_j(\rho)\,,
\end{align}
and differentiating \(Q_j\) yields
\(Q_j'(t) = \left(j-\tfrac12\right)t^{j-1}g(t) - jt^jg(t)-\tfrac12t^{j-1}\).
Integrating this identity over \([0,h]\) and using \eqref{eq:even-readable-integral} gives
$$
   h^{j-1/2}|\rho|\arcsin|\rho|
	   = Q_j(h)-Q_j(0)
     = \left(j-\tfrac12\right)h^{j-1}M_{j-1}(\rho) - jh^jM_j(\rho)-\tfrac{h^j}{2j}.
$$
The recursion in \eqref{eq:even-M-recursion} then follows upon division by \(jh^j\).

We now rewrite \(\mathcal A_{2k}\) in \eqref{thm:general:eq} in terms of \((M_j)_{j \in \{0,k,\dots,2k-1\}}\) and a remainder term.
By \citep[15.4.4 \& 16.12.2]{DLMF},
$$
  (\arcsin\sqrt t)^2
    = \sum_{n=1}^{\infty} r_n \, t^n,
		\qquad \qquad 
		r_n = \frac{4^n}{2n^2\binom{2n}{n}}\,.
$$
Using that \(g(t)=\frac{\pi}{2 \, \sqrt{t(1-t)}} - \tfrac{\dd}{\dd t} (\arcsin\sqrt t)^2\) together with \eqref{eq:even-readable-integral2} gives
$$
   M_j(\rho)
     = \frac{\pi}{2}h^{-j} \int_{[0,h]} t^{j-1/2}(1-t)^{-1/2} \dd \lambda(t)
	     - h^{-j} \int_{[0,h]} t^j \, \frac{\dd}{\dd t} (\arcsin\sqrt t)^2 \dd \lambda(t).
$$
Combining the binomial series
\((1-t)^{-1/2}
  = \sum_{\ell=0}^{\infty}\frac{(1/2)_\ell}{\ell!}t^\ell\) with 
\(\tfrac{\dd}{\dd t} (\arcsin\sqrt t)^2
  = \sum_{n=1}^{\infty} n \, r_n \, t^{n-1}\),
termwise integration yields
\begin{align} \label{eq:even-readable-series}
  M_j(\rho)
  & = \frac{\pi}{2}h^{-j} \sum_{\ell=0}^{\infty} \frac{(1/2)_\ell}{\ell!} \frac{h^{j+\ell+1/2}}{j+\ell+1/2}
	    - h^{-j} \sum_{n=1}^{\infty} n \, r_n \, \frac{h^{j+n}}{j+n} \notag
	\\
  & = \frac{\pi}{2} \sum_{\ell=0}^{\infty} \frac{(1/2)_\ell}{\ell!(j+\ell+1/2)}h^{\ell+1/2}
      - \sum_{n=1}^{\infty} \frac{n \, r_n}{n+j} \, h^n.
\end{align}
The interchange of summation and integration is justified by the uniform convergence of both series for \(0<h<1\).
To derive the announced decomposition of \(\mathcal A_{2k}\), define
$$
  R_k(z)
  = \frac{(z-2k)_k}{(z+k)_k}.
$$
Since \((z+k)_k=\prod_{r=k}^{2k-1}(z+r)\), the function \((-1)^k \, R_k(z)/z\) has simple poles at \(z \in \{0,-k,\ldots,-(2k-1)\}\). Its residue at \(z=0\) is
\((-1)^k \, R_k(0)
  = (-1)^k \frac{(-2k)_k}{(k)_k} = 2\).
For \(k\leq j\leq2k-1\), the residue at \(z=-j\) is given by
$$
  (-1)^k \, \frac{(-j-2k)_k}{(-j)\prod_{\substack{r=k\\r\neq j}}^{2k-1}(r-j)}
	  = (-1)^{k+j+1} \frac{(j+2k)!}{j(j+k)!(j-k)!(2k-1-j)!}
	  = c_{k,j}.
$$
Therefore,
\begin{align}\label{eq:even-readable-pf}
 (-1)^k\frac{R_k(z)}z
 &=\frac2z+\sum_{j=k}^{2k-1}\frac{c_{k,j}}{z+j}.
\end{align}
Using the identities
\(\frac{(k+1/2)_\ell(1/2-k)_\ell}
      {(2k+1/2)_\ell(1/2-2k)_\ell}
  = (-1)^k \, R_k(\ell+1/2)\) and 
\(\frac{(1/2)_\ell}{(3/2)_\ell}
  = \frac1{2\ell+1}\),
the series \(\mathcal A_{2k}\) becomes
\begin{align} \label{eq:even-readable-pf2}
  \pi\sqrt h\,\mathcal A_{2k}(h)
  & = \pi \sqrt h \,{}_4F_3 (\tfrac12,\tfrac12,\tfrac{2k+1}{2},\tfrac{1-2k}{2}; \tfrac32,2k+\tfrac12,\tfrac12-2k; h) \notag
	  = \pi \sqrt h \, \sum_{\ell=0}^{\infty} \frac{(\frac12)_\ell^2}{(\frac32)_\ell\,\ell!} \, \frac{(k+\frac12)_\ell(\frac12-k)_\ell}
{(2k+\tfrac12)_\ell(\tfrac12-2k)_\ell} h^\ell \notag
  \\
  & = \frac{\pi}{2}\sum_{\ell=0}^{\infty} \frac{(\tfrac12)_\ell}{\ell!}h^{\ell+1/2} \frac{(-1)^k \, R_k(\ell+\tfrac12)}{\ell+\tfrac12} \notag
	  = \frac{\pi}{2}\sum_{\ell=0}^{\infty} \frac{(\tfrac12)_\ell}{\ell!}h^{\ell+1/2} \left( \frac2{\ell+1/2} + \sum_{j=k}^{2k-1}
       \frac{c_{k,j}}{\ell+j+1/2} \right) \notag
	\\
	& = 2M_0(\rho) + \sum_{j=k}^{2k-1} c_{k,j} M_j(\rho) 
	    + \sum_{n=1}^{\infty} \left( 2 \, r_n \, h^n + \sum_{j=k}^{2k-1} \frac{n \, r_n \, c_{k,j}}{n+j} \, h^n \right) \notag
	\\
	& = 2M_0(\rho) + \sum_{j=k}^{2k-1} c_{k,j} M_j(\rho) 
	    + (-1)^k \, \sum_{n=1}^{\infty} r_n \, R_k(n) \, h^n\,,
\end{align}
where the third last identity follows from \eqref{eq:even-readable-pf}, the second last from \eqref{eq:even-readable-series}, and the last identity from applying \eqref{eq:even-readable-pf} at \(z=n\).

We now express the remainder term in \eqref{eq:even-readable-pf2} in terms of \(\mathcal B_{2k}\) by matching the coefficients of the corresponding series.
For \(1\leq n \leq k\), canceling common factors in \(R_k(n)\) yields
$$
 (-1)^k \, r_n \, R_k(n)
   = r_n \, \frac{\binom{2k-n}{k}}{\binom{2k+n-1}{k}}
   = b_{k,n}.
$$
Since \(R_k(n)=0\) for \(k+1\leq n \leq 2k\), it remains to identify the tail beginning at \(n=2k+1\).
Define \(t_n = r_n \, R_k(n)\).
Taking ratios of consecutive terms gives
$$
 \frac{t_{n+1}}{t_n}
   = \frac{n^2}{(n+1)(n+1/2)} \, \frac{(n-k)(n+k)}{(n-2k)(n+2k)}
   \qquad \textrm{ and } \qquad t_{2k+1}
   = \frac{2^{4k+1}k!(3k)!((2k)!)^2}{(4k+2)!(4k)!}.
$$
Replacing \(n=2k+1+\ell\), the ratio becomes
$$
 \frac{(2k+1+\ell)^2(k+1+\ell)(3k+1+\ell)}{(2k+2+\ell)(2k+3/2+\ell)(\ell+1)(4k+1+\ell)}\,.
$$
Next, consider the series expansion \(\mathcal B_{2k}(h) = \sum_{\ell=0}^{\infty} a_\ell \, h^\ell\), whose coefficients are, by definition,
$$
  a_\ell
	= \frac{(2k+1)_\ell^2(k+1)_\ell(3k+1)_\ell}{(2k+3/2)_\ell(2k+2)_\ell(4k+1)_\ell\,\ell!}.
$$
Applying the identity \((k)_{\ell+1}=(k+\ell)(k)_\ell\) then yields
\( a_{\ell+1} / a_\ell = t_{2k+2+\ell} / t_{2k+1+\ell}\).
Since \(a_0=1\), it follows that \(t_{2k+1+\ell}=t_{2k+1}a_\ell\). 
Consequently, 
$$
   \sum_{n=2k+1}^{\infty} t_nh^n
     = t_{2k+1}h^{2k+1} \sum_{\ell=0}^{\infty}a_\ell h^\ell
	   = t_{2k+1}h^{2k+1} \, \mathcal B_{2k}(h).
$$
The gamma-function identities in the definition of \(e_\nu\) further give
\(e_{2k} = -\frac{3(-1)^k}{\pi^2}\,t_{2k+1}\).
Since \(\mathcal B_{2k}(0)=1\), combining these identities yields
$$
  (-1)^k \sum_{n=2k+1}^{\infty} r_n \, R_k(n) \, h^n
  = -\frac{\pi^2}{3}e_{2k}h^{2k+1}\mathcal B_{2k}(h)\,,
$$
so that \eqref{eq:even-readable-pf2} becomes
$$
  2M_0(\rho)+\sum_{j=k}^{2k-1}c_{k,j}M_j(\rho) + \sum_{\ell=1}^{k}b_{k,\ell}h^\ell
  = \pi\sqrt h\,\mathcal A_{2k}(h)
    + \frac{\pi^2}{3}e_{2k}h^{2k+1}\mathcal B_{2k}(h).
$$
By Theorem~\ref{thm:general}, the right-hand side equals \(\pi^2(1-\xi(2k,\rho))/3\), which proves \eqref{eq:xi-even}.

Finally, continuity of the functions \(M_j\) on \((-1,1)\) together with Proposition \ref{prop.xi.cont}, extends the identity to \(\rho=0\) (\(h=1\)).
The continuous extension to \(\rho=\pm1\) (\(h=0\)), where \(\xi(2k,\pm1)=1\), follows directly from Theorem~\ref{thm:general}.
\end{proof}

By Lemma~\ref{eq:xi-expectation},
\begin{align} \label{Thm:Monot.rho.1}
  \xi(\nu,\rho)
  & = 1-\frac{3\sqrt h}{\pi} \int_{[0,1]} \E[(1+hu^2W_m)^{m-1}] \dd \lambda(u)
	  = 1-\frac{3\sqrt h}{\pi} \int_{[0,1]} \Phi_m(hu^2) \dd \lambda(u)\,,
\end{align}
where \(\Phi_m: [0,1] \to \R\) is defined by 
\begin{align} \label{Thm:Monot.rho.2}
  \Phi_m(s):=\mathbb E(1+sW_m)^{m-1}\,,
\end{align}
and where \(W_m=B_m D_m\) with independent \(B_m \sim\beta\left(m,m-\frac12\right)\) and \(D_m \sim \beta'\left(\frac12,2m-\frac12\right)\).
The properties of \(\Phi_m\) play a central role in the proofs of Theorem~\ref{Thm:MonCCX} and Theorem~\ref{Thm:Monot.Nu} below.

\begin{lemma}\label{eq:xi.Phi}
For every \(m=(\nu+1)/2>1/2\), the function $(m,s) \mapsto \Phi_m(s)$ is finite, continuous, and strictly positive on $(1/2,\infty) \times [0,1]$.
Moreover, for every \(s \in (0,1]\), the mapping \(m \mapsto \Phi_m(s)\) is strictly increasing on \((1/2,\infty)\).
\end{lemma}

\begin{proof}
By \eqref{eq:W_moments}, $\Phi_m$ is finite and strictly positive since $W_m>0$ almost surely. Continuity follows from dominated convergence.

We prove the strict monotonicity by comparing the random variables \(W_m\) and weighted versions of these under the usual stochastic order.  In particular, we show that \(W_n\leq_{\mathrm{st}}W_m\) for \(n>m\), where $X_1\leq_{\mathrm{st}}X_2$ means $P(X_1>t) \leq P(X_2>t)$ for all $t\in\R$. 
To establish this, we show that the density ratio \(f_{W_n}/f_{W_m}\) is unimodal for \(n > m\), with \(\lim_{w\downarrow0} f_{W_n}(w)/f_{W_m}(w) > 1\) and \(\lim_{w \to \infty} f_{W_n}(w)/f_{W_m}(w) = 0\). 
Hence, the ratio crosses the level \(1\) at most once, necessarily from above to below. The stochastic ordering then follows from the single-crossing criterion in \cite[Proposition~2.1]{KlenkeMattner2010}.

To establish unimodality, we first decompose the density \(f_{W_{m}}\).  
By \cite[Corollary 3.2.1]{NagarZarrazola2004},
\begin{align}\label{Thm:Monot.Nu:Jm}
  f_{W_{m}}(w)
	& = \underbrace{\frac{1}{\B(m,m-\frac12) \B(\frac12,2m-\frac12)}}_{=: K_m} \, \frac{1}{\sqrt{w}} \, \underbrace{\int_{[0,1]} v^{3m-3/2}(1-v)^{m-3/2} \, \left(1+\frac vw\right)^{-2m} \, \frac{1}{w^{2m}} \dd \lambda(v)}_{=: P_m(w) >0}.	
\end{align}
The function \(w \mapsto P_m(w)\) is completely monotone and hence log-convex by \cite[Lemma 4.3]{Merkle2014}.
Consequently, 
$ \eta_m(w)
  := - P_m'(w) / P_m(w) $ is nonincreasing on $(0,\infty)$.  Moreover, \(0<\eta_m(w)\le\frac{2m}{w}\) so that \(\lim_{w\to\infty}\eta_m(w)=0\).
By the integral representation \cite[Eq.~15.1.2 \& 15.6.1]{DLMF}, \(P_m\) satisfies
$$
  P_m(w)
  = B\left(3m-\tfrac12,m-\tfrac12\right) \, w^{-2m} \, {}_2F_1\left(2m,3m-\tfrac12;4m-1;-\tfrac1w \right).
$$
Consequently, the hypergeometric differential equation for \({}_2F_1\) \cite[Eq.~15.10.1]{DLMF}, together with the change of variables $z=-1/w$, yields
\begin{align}\label{eq:df-ode}
  0
  & = w(1+w) \, P_m''(w) + \left(3w+\tfrac32-m\right) \, P_m'(w) - 4m(m-1) \, P_m(w).
\end{align}
We next use \eqref{eq:df-ode} to compare the functions \(P_m\) for different parameter values, thereby drawing on the properties of \(\eta_m\) established above.
Let $n>m$ and define
$$
  R(w)
	  = \frac{P_n(w)}{P_m(w)}
	\qquad \textrm {and } \qquad  
	p(w) 
	  = w^{3/2-n}(1+w)^{3/2+n}P_m(w)^2. 
$$
Substituting \(P_n=RP_m\) into \eqref{eq:df-ode} and using the corresponding equation for \(P_m\) gives
\begin{align} \label{eq:df-ode2}
  0
	& = R''(w) + \left(2\frac{P_m'(w)}{P_m(w)}+\frac{3w+\frac32-n}{w(1+w)}\right) R'(w) + \frac{n-m}{w(1+w)} \left[\eta_m(w)-4(n+m-1)\right] R(w).
\end{align}
Since $p(w)>0$ for every $w>0$ and
\begin{align*}
  \frac{p'(w)}{p(w)}
  & = 2 \, \frac{P_m'(w)}{P_m(w)}+\frac{\frac32-n}{w} + \frac{\frac32+n}{1+w}
	  = 2 \, \frac{P_m'(w)}{P_m(w)}+\frac{3w+\frac32-n}{w(1+w)}\,,
\end{align*}
\eqref{eq:df-ode2} then yields 
$$
  (p(w)R'(w))'
	  = p(w) R''(w) + p'(w)R'(w)
	  = -\frac{(n-m)p(w)R(w)}{w(1+w)} \left[\eta_m(w)-4(n+m-1)\right].
$$
Since $n>m$, $p(w)>0$, and $R(w)>0$, the sign of $(pR')'(w)$ is determined entirely by the expression \(\eta_m(w)-4(n+m-1)\). 
As $\eta_m$ is nonincreasing and converges to zero, this expression can change sign at most once, from positive to negative.
Consequently, $pR'$ is initially nonincreasing and then nondecreasing, with the initial part possibly absent. 
This observation is key to establishing the unimodality of the density ratio.
To determine the sign of $pR'$ on its final nondecreasing segment, we examine the behavior of $R(w)$ as $w\to\infty$.
For each fixed $v\in(0,1)$, we have \(\left(1+\frac vw\right)^{-2m} \to 1\) as \(w\to\infty\).
Thus, dominated convergence in \eqref{Thm:Monot.Nu:Jm} yields
\( \lim_{w\to\infty} w^{2m}P_m(w) = \B \left(3m-\frac12,m-\frac12\right)\).
Since $n>m$, it follows that, as \(w\to\infty\), 
$$
  R(w)
    = \frac{P_n(w)}{P_m(w)} 
	  \approx \frac{\B(3n-\frac12,n-\frac12)}{\B(3m-\frac12,m-\frac12)}w^{-2(n-m)}
		\to 0.
$$
Therefore, $pR'$ must be negative on its final nondecreasing segment.  Otherwise, since $p(w)>0$, $R$ would eventually increase, contradicting $R(w)>0$ and $R(w)\to0$. 
Thus, $R$ is either nonincreasing throughout $(0,\infty)$ or first nondecreasing and then nonincreasing, and hence unimodal in either case.  Since the density ratio $f_{W_n}/f_{W_m}$ is a positive multiple of $R$, it is also unimodal.

We now examine the behavior of the density ratio \(f_{W_n}/f_{W_m}\) at both endpoints.
For $w\downarrow0$, the integrand defining $P_m(w)$ in \eqref{Thm:Monot.Nu:Jm} converges pointwise to 
\(v^{m-3/2}(1-v)^{m-3/2}\).
Since $m>\frac12$, this function is integrable on $(0,1)$ and dominates the original integrand.
Dominated convergence in \eqref{Thm:Monot.Nu:Jm} therefore gives
$$
  \lim_{w\downarrow0}P_m(w)
    = \int_{[0,1]} v^{m-3/2}(1-v)^{m-3/2} \dd \lambda(v) 
    = \B\left(m-\tfrac12,m-\tfrac12\right).
$$
Consequently, as \(w\downarrow0\),
\(f_{W_m}(w) \approx K_m \, w^{-1/2} \, \B\left(m-\tfrac12,m-\tfrac12\right) = 2c_m \, w^{-1/2}\),
with \(c_m\) as defined in \eqref{Abbrev.Short}.
Since \(c_m\) is strictly increasing as a function of $m$, for $n>m$, 
\(\lim_{w\downarrow0} f_{W_n}(w) / f_{W_m}(w) = c_n / c_m > 1\). 
On the other hand, since $K_n/K_m>0$ and $R(w)\to0$,
\(\lim_{w\to\infty} f_{W_n}(w) / f_{W_m}(w) = \lim_{w\to\infty} (K_n / K_m) \, R(w) = 0\).
We have therefore shown that \(W_n\leq_{\mathrm{st}}W_m\) for \(n>m\).

We next extend the stochastic comparison to weighted versions of \(W_m\).
For $k \in \N_0$ and $m>\max\{k,\frac12\}$, let $W_{m,k}$ have density \(f_{W_{m,k}}: (0,\infty) \to \R\) given by
$$
  f_{W_{m,k}}(w)
	= \frac{w^k\,f_{W_m}(w)}{\E[W_m^k]}.
$$
By \eqref{eq:W_moments}, $\E[W_m^k]<\infty$ because $k<2m-\frac12$.
The order just proved implies $\E[W_n^k] \leq \E[W_m^k]$ for all moments considered below.
Moreover,
$$
  \frac{f_{W_{n,k}}(w)}{f_{W_{m,k}}(w)}
  = \frac{\E[W_m^k]}{\E[W_n^k]} \, \frac{f_{W_n}(w)}{f_{W_m}(w)}.
$$
This ratio retains the unimodal shape of \(f_{W_n} / f_{W_m}\) and satisfies
\(\lim_{w\downarrow0} f_{W_{n,k}}(w) / f_{W_{m,k}}(w) > 1\) and \linebreak
\(\lim_{w\to\infty} f_{W_{n,k}}(w) / f_{W_{m,k}}(w) = 0\).
Applying \cite[Proposition~2.1]{KlenkeMattner2010} again yields \(W_{n,k} \leq_{\mathrm{st}} W_{m,k}\) whenever \(\max\{k,1/2\}<m<n\). Thus, \(\E[\phi(W_{n,k})] \leq \E[\phi(W_{m,k})]\) for every nondecreasing function \(\phi\) for which the expectations exist.

To establish the claimed monotonicity of \(m \mapsto \Phi_m\), we analyze the derivative of \(\Phi_m\).
For $k \in \N_0$ and $m>\max\{k,\frac12\}$, differentiation under the expectation is justified by finite moments and yields
$$
  \Phi_m^{(k)}(s)
  = \E \left[ W_m^k(1+sW_m)^{m-1-k} \right] \, \prod_{\ell=1}^k(m-\ell).
$$
Since the moments are finite, these derivatives also exist at $s=0$.
Define \(A_0(m) = 1\) and, for \(k\in \N\), \(A_k(m) = \E[W_m^k] \, \prod_{\ell=1}^k(m-\ell)\).
Then
\begin{align} \label{eq:phi-derivatives}
  \Phi_m^{(k)}(s)
  & = A_k(m) \, \E[(1+sW_{m,k})^{m-1-k}].
\end{align}
In particular, $\Phi_m^{(k)}(0)=A_k(m)$. 
We compare the two factors separately. First we establish monotonicity of $m \mapsto A_k(m)$, and then we use the stochastic comparison of $W_{m,k}$ to compare the corresponding expectations in \eqref{eq:phi-derivatives}.
By \eqref{eq:W_moments}, the coefficients $A_k(m)$ satisfy
\begin{align}\label{eq:df-coefficients}
  A_k(m)
  & = \left(\tfrac12\right)_k \, \prod_{j=0}^{k-1} \frac{m+j}{4m(m-1)+\frac34-j(j+1)} \, \prod_{\ell=1}^k \, (m-\ell)
	  = \left(\frac12\right)_k \prod_{j=0}^{k-1} \frac{\lambda_m-j(j+1)}{4\lambda_m+\frac34-j(j+1)},
\end{align}
where we define \(\lambda_m = m(m-1)\).
The advantage of this representation is that the dependence on $m$ now occurs only through $\lambda_m$.
Every factor in \eqref{eq:df-coefficients} is strictly positive and strictly increasing in $\lambda_m$.
Since $(\frac12)_k>0$, it also follows that the product in \eqref{eq:df-coefficients} is strictly increasing in $\lambda_m$.
Since $m \mapsto \lambda_m$ itself is strictly increasing, we conclude that \(m \mapsto A_k(m)\) is strictly increasing on \((k,\infty)\) for every $k \in \N$.
Now fix $k \in \N_0$ and \(\max\left\{k,\frac12\right\}<m<n<k+1\).
Since \(m-1-k, n-1-k \in (-1,0)\), \(W_{n,k} \leq_{\mathrm{st}} W_{m,k}\) implies
$$
  \E[(1+sW_{n,k})^{n-1-k}]
   \geq \E[(1+sW_{m,k})^{n-1-k}]
    >   \E[(1+sW_{m,k})^{m-1-k}]
$$
for every $s > 0$.
The second inequality is strict because $W_{m,k}>0$ almost surely and $n>m$.
Since \(A_k(n)\ge A_k(m)>0\) for every \(k\in\mathbb N_0\), with equality when \(k=0\), the preceding strict inequality yields \(\Phi_n^{(k)}(s)>\Phi_m^{(k)}(s)\) for every \(s>0\).

For $k=0$, this already gives $\Phi_m(s) < \Phi_n(s)$ for all \(1/2 < m < n < 1\). 
For $k \in \N$, the derivatives in \eqref{eq:phi-derivatives} satisfy
$$
  \Phi_n^{(j)}(0)
      =  A_j(n)
	  \geq A_j(m)
      = \Phi_m^{(j)}(0)
$$
for all \(j \in \{0,\ldots,k-1\}\).
Define \(\delta:=\Phi_n-\Phi_m\).  
Then \(\delta(0)=0\) while \(\delta^{(j)}(0) \geq 0\) for all \(j \in \{1,\ldots,k-1\}\).
Moreover, the strict comparison of the \(k\)-th derivatives established above yields 
\(\delta^{(k)}(s) > 0\) for every \(s>0\).
Taylor's formula with integral remainder gives, for $s>0$,
$$
  \delta(s)
  = \sum_{j=0}^{k-1} \frac{\delta^{(j)}(0)}{j!}s^j + \frac{1}{(k-1)!} \int_{[0,s]} (s-t)^{k-1}\delta^{(k)}(t) \dd \lambda(t).
$$
Consequently, \(\delta(s)>0\) for all \(s>0\), which proves \(\Phi_m(s)<\Phi_n(s)\) for all \(s>0\).

Continuity extends monotonicity across the endpoints \(\max\left\{k,\frac12\right\}<m<n<k+1\), \(k \in \N_0\), which implies that $m \mapsto \Phi_m(s)$ is nondecreasing on the whole of $(1/2,\infty)$ for every $s>0$.
Since every interval with distinct endpoints in $(1/2,\infty)$ contains a nonempty subinterval of one of the above open intervals, the increase is strict on the whole of $(1/2,\infty)$.
\end{proof}

\begin{proof}[\textbf{Proof of Theorem~\ref{Thm:MonCCX}}] 
Consider \(r=|\rho|\). The substitution $t=\sqrt{h}\,u=\sqrt{1-r^2}\,u$ in \eqref{Thm:Monot.rho.1} gives
$$
  \xi(\nu,r)
  = 1-\frac{3}{\pi} \int_{[0,\sqrt{1-r^2}]} \Phi_m (t^2) \dd \lambda(t)\,.
$$
Hence, for $0\le r_1<r_2<1$, \(\xi(\nu,r_2)-\xi(\nu,r_1) >0\),
since $\Phi_m$ is strictly positive by Lemma~\ref{eq:xi.Phi}.
Finally, the integral representation shows that \(\xi(\nu,r)<1\) for every \(r<1\), since \(\Phi_m\) is strictly positive. As \(\xi(\nu,1)=1\), strict monotonicity extends to \(r=1\).
\end{proof}

\begin{proof}[\textbf{Proof of Theorem~\ref{Thm:Monot.Nu}}] 
The result is immediate from the representation of \(\xi\) in Lemma~\ref{eq:xi-expectation} and the fact that \(m \mapsto \Phi_m(s)\) is strictly increasing for every \(s>0\), as established in Lemma~\ref{eq:xi.Phi} 
\end{proof}

\begin{proof}[\textbf{Proof of Theorem \ref{Thm:NuToInfty}}]
Since $T_\nu\to\Phi$ as $\nu\to\infty$ and $\Phi$ is continuous and strictly increasing, we have
\(q_\nu(u) \to \Phi^{-1}(u)\) for all \(u\in(0,1)\).
Consequently, as $\nu\to\infty$,
$$
  s_{\nu,\rho}(q_\nu(u))
  = \sqrt{(1-\rho^2) \frac{\nu+q_\nu(u)^2}{\nu+1}} \to \sqrt{1-\rho^2}.
$$
By Proposition~\ref{prop:conditional}, for every
$(u,v)\in(0,1)^2$,
$$
  \partial_1 C_{\nu,\rho}(u,v)
    = T_{\nu+1}\left(\frac{q_\nu(v)-\rho q_\nu(u)}{s_{\nu,\rho}(q_\nu(u))} \right)
    \to \Phi\left( \frac{\Phi^{-1}(v)-\rho\Phi^{-1}(u)}{\sqrt{1-\rho^2}}\right)
    = \partial_1 C_\rho(u,v)
$$
and, since $0\leq\partial_1 C_{\nu,\rho}\leq1$, also in \(L^p\).
The result for \(\xi\) now follows from \cite[Theorem 3.5]{AnsariFuchs2026}.
\end{proof}

\begin{proof}[\textbf{Proof of Theorem \ref{Thm:NuTo0}}] 
We first show that, for $1/2<a<b<1$, \(\frac{q_\nu(a)}{q_\nu(b)}\to 0\).  The proof relies on \cite[Appendix C]{McNeil2022} which uses the stochastic decomposition of Student-\(t\) distributions into radial and spherical components.  Following their notation, let $G_{2,\nu}$ denote the distribution function of the radial component of a bivariate Student-$t$ distribution.
Then, by \cite[Proposition~8]{McNeil2022}, for any $\varepsilon>0$,
$ \lim_{\nu\downarrow0} G_{2,\nu}\bigl(\varepsilon q_\nu(b)\bigr)
  = 2b-1>2a-1
  = \lim_{\nu\downarrow0} G_{2,\nu}\bigl(q_\nu(a)\bigr).$
Monotonicity of $G_{2,\nu}$ thus gives $q_\nu(a) < \varepsilon q_\nu(b)$ for all sufficiently small $\nu$, hence \(\frac{q_\nu(a)}{q_\nu(b)}\to 0\).
Define $s(u)=\min\{u,1-u\}$.
By symmetry, $q_\nu(1-t)=-q_\nu(t)$ and hence $|q_\nu(t)|=q_\nu(1-s(t))$ for all $t\ne1/2$.
Consequently,
$$
  \frac{|q_\nu(v)|}{|q_\nu(u)|}
  = \frac{q_\nu(1-s(v))}{q_\nu(1-s(u))} 
	\to \begin{cases}
      0 & s(v)>s(u)
			\\
      +\infty 
			& s(v)<s(u).
\end{cases}
$$
Taking into account the sign of $q_\nu(v)$ yields, for almost every $(u,v)\in(0,1)^2$,
\begin{align*}
\frac{q_\nu(v)}{|q_\nu(u)|}\longrightarrow
\begin{cases}
  -\infty
	& v<s(u)
	\\
  0 
	& s(u)<v<1-s(u)
	\\
  +\infty 
	& v>1-s(u).
\end{cases}
\end{align*}
Moreover, \cite[Lemma~3]{McNeil2022} implies $|q_\nu(u)|\to\infty$ for $u\ne1/2$,
and hence, as \(\nu\downarrow 0\), 
\begin{align*}
  \frac{s_{\nu,\rho}(q_\nu(u))}{|q_\nu(u)|}
	& = \sqrt{\frac{1-\rho^2}{\nu+1} \left(1+\frac{\nu}{q_\nu(u)^2}\right)}
  \longrightarrow\sqrt{1-\rho^2}.
\end{align*}
Now, using Proposition~\ref{prop:conditional} and $T_{\nu+1}\to T_1$, we obtain
\begin{align*}
  \lim_{\nu\downarrow0}\partial_1 C_{\nu,\rho}(u,v)
  & = \lim_{\nu\downarrow0} T_{\nu+1}\left( \frac{\displaystyle\frac{q_\nu(v)}{|q_\nu(u)|}-\rho\frac{q_\nu(u)}{|q_\nu(u)|}}{\displaystyle\frac{s_{\nu,\rho}(q_\nu(u))}{|q_\nu(u)|}} \right)
	  = \begin{cases}
    0 
		& v<s(u)
		\\
    \displaystyle \frac12- \frac{\operatorname{sgn}(u-1/2)}{\pi}\arcsin\rho
    & s(u)<v<1-s(u)
		\\
    1
		& v>1-s(u).
  \end{cases}
	\\
  & = \alpha_\rho\ind_{\{u<v\}} + (1-\alpha_\rho)\ind_{\{u+v>1\}}
	  = \partial_1 C_{0,\rho}(u,v)\,,
\end{align*}
where the second last identity follows by a distinction of cases.
Thus,
\(\partial_1 C_{\nu,\rho} \to \partial_1 C_{0,\rho}\)
almost everywhere on \((0,1)^2\) and, since $0\leq\partial_1 C_{\nu,\rho}\leq1$, also in \(L^p\).
The result for \(\xi\) now follows from \cite[Theorem 3.5]{AnsariFuchs2026}.
\end{proof}

\section{Proofs of Section \ref{sec:estimator}} \label{sec:Proof3}

\begin{proof}[\textbf{Proof of Lemma~\ref{lem:joint-criterion}}]
Identifiability of the family \(\{C_{\nu,\rho}:(\nu,\rho)\in\Theta\}\) follows from Corollary \ref{cor:df-identifiability}, while uniform convergence follows from \cite[Theorem~1(a) \& Corollary 1]{Liu2026}.
Finally,
$$
  L(\theta_0)-L(\theta)
  =
  \int_{(0,1)^2}
  c_{\theta_0}(u,v)
  \log\frac{c_{\theta_0}(u,v)}{c_\theta(u,v)}
  \dd\lambda(u,v)
  \geq 0
$$
by the Kullback--Leibler inequality, with equality if and only if \(c_\theta=c_{\theta_0}\) almost everywhere. 
By identifiability in Corollary~\ref{cor:df-identifiability}, this equality implies \(\theta=\theta_0\).  Thus, \(L\)  is uniquely maximized at \(\theta_0\).
\end{proof}

\begin{proof}[\textbf{Proof of Theorem~\ref{thm:joint-consistency}}] 
The strong law of large numbers for bounded U-statistics gives \(\widehat\tau_n \xrightarrow{\mathrm{a.s.}}\tau_0\).
Since \(\rho=\sin(\pi\tau / 2)\), continuity yields \(\sin(\pi\widehat\tau_n/2) \xrightarrow{\mathrm{a.s.}}\rho_0\).
Because \(|\rho_0|<r\), truncation to \([-r,r]\) is almost surely inactive for sufficiently large \(n\).
Consequently, \(\widehat\rho_n^{\mathrm{TP}} \xrightarrow{\mathrm{a.s.}}\rho_0\).
By Lemma~\ref{lem:joint-criterion} and uniform continuity of \(L\) on the
compact set \(\Theta\),
\begin{align*}
  \sup_{\nu\in[a,b]}
  \left| L_n(\nu,\widehat\rho_n^{\mathrm{TP}}) - L(\nu,\rho_0) \right|
  & \leq \sup_{\theta\in\Theta} |L_n(\theta)-L(\theta)| 
	       + \sup_{\nu\in[a,b]} \left|L(\nu,\widehat\rho_n^{\mathrm{TP}}) - L(\nu,\rho_0) \right|
    \xrightarrow{\mathrm{a.s.}}0.
\end{align*}
The function \(\nu \mapsto L(\nu,\rho_0)\) is continuous on \([a,b]\) and, by Lemma~\ref{lem:joint-criterion}, has its unique maximum at \(\nu_0\).
Hence \cite[Proposition~5.1]{Liebscher2009} yields \(\widehat\nu_n^{\mathrm{TP}}\xrightarrow{\mathrm{a.s.}}\nu_0\).
Thus \(\widehat\theta_n^{\mathrm{TP}} \xrightarrow{\mathrm{a.s.}}\theta_0\).
The claimed strong consistency now follows from the continuity of \(\xi\), as shown in Proposition~\ref{prop.xi.cont}.
\end{proof}

\begin{proof}[\textbf{Proof of Lemma~\ref{lem:joint-target-derivative}}]
Let $m=(\nu+1)/2$ and $h=1-\rho^2$ as defined in \eqref{Abbrev.Short}.
By Lemma~\ref{eq:xi-expectation}, for $0<h\le1$,
$$
  \xi(\nu,\rho)
	= 1 - \frac{3\sqrt h}{\pi} \, \Psi(m,h),
	\quad \textrm{ where } \quad 	
	\Psi(m,h)
	:= \int_{(0,1)} \int_{(0,\infty)} (1+hu^2w)^{m-1} \, f_{W_m}(w) \dd \lambda(w) \dd \lambda(u)\,,
$$
with \(f_{W_m}(w)=K_mw^{-1/2}P_m(w)\) as given in \eqref{Thm:Monot.Nu:Jm}.
We extend \(\Psi\) to \(h>0\) and show that it is continuously differentiable on \((1/2,\infty) \times (0,\infty)\), so that \(h=1\) is an interior point of its domain.
The assertion then follows by applying the chain rule to
$$ \xi(\nu,\rho)
   = 1-\frac{3\sqrt{1-\rho^2}}{\pi} \Psi\left(\frac{\nu+1}{2},1-\rho^2\right).$$

Fix a compact rectangle $[m_0,m_1] \times [h_0,h_1] \subseteq (1/2,\infty) \times (0,\infty)$.
Differentiating the integrand $p_m$ of $P_m$ with respect to $m$ yields
\(\partial_m p_m(v,w) = (3\log v+\log(1-v)-2\log(v+w)) \, p_m(v,w)\).
For $0<w<1$, since $v\le v+w\le2$, 
$$
  p_m(v,w) + |\partial_m p_m(v,w)|
  \le C \, v^{m_0-3/2}(1-v)^{m_0-3/2} \, \big(1+|\log v|+|\log(1-v)|\big).
$$
For $w\ge1$, using $v+w\ge w$ and $\log(v+w)\le\log(2w)$ yields
$$
  p_m(v,w)+|\partial_mp_m(v,w)|
  \le C \, w^{-2m} \, v^{3m_0-3/2}(1-v)^{m_0-3/2} \, \big(1 +|\log v|+|\log(1-v)|+\log w \bigr).
$$
Both bounds are integrable with respect to $v$ since $m_0>1/2$.
Thus, differentiation under the integral defining $P_m$ is justified and yields
\begin{align*}
  0 < P_m(w)
	& \le K
    \begin{cases}
      1 & 0<w<1\\
      w^{-2m} & 1\le w < \infty
    \end{cases}\\
    |\partial_mP_m(w)|&\le K
    \begin{cases}
      1 & 0<w<1\\
      w^{-2m}(1+\log w) & 1\le w < \infty
    \end{cases}
\end{align*}
for some \(K < \infty\).
The same majorant also ensures joint continuity of $P_m(w)$ and $\partial_mP_m(w)$ on every compact rectangle contained in \((1/2,\infty)\times(0,\infty)\), and hence throughout the entire domain. 
Since $K_m$ in \eqref{Thm:Monot.Nu:Jm} is continuously differentiable for $m>1/2$, both $K_m$ and $K_m'$ are bounded on $[m_0,m_1]$, hence 
$\partial_m f_{W_m}(w)
	= w^{-1/2} \bigl( K_m' P_m(w) + K_m \partial_m P_m(w) \bigr)$.
Therefore,
\begin{align*}
  0 < f_{W_m}(w) 
	& \le K
    \begin{cases}
    w^{-1/2} & 0<w<1\\
    w^{-2m-1/2} & 1\le w < \infty
    \end{cases}\\
    |\partial_mf_{W_m}(w)|&\le K
    \begin{cases}
    w^{-1/2} & 0<w<1\\
    w^{-2m-1/2}(1+\log w) &1\le w < \infty
    \end{cases}
\end{align*}
for some \(K < \infty\).

Denote by \(\psi_{m,h}\) the integrand of \(\Psi\).
Then
\begin{align*}
  \partial_m \psi_{m,h}(u,w)
  & = (1+hu^2w)^{m-1} \bigl[ \partial_mf_{W_m}(w) + f_{W_m}(w)\log(1+hu^2w) \bigr],
	\\
  \partial_h \psi_{m,h}(u,w) 
  & = (m-1) \, \frac{u^2w}{1+hu^2w} \, \psi_{m,h}(u,w).
\end{align*}
On the compact rectangle \([m_0,m_1] \times [h_0,h_1]\), 
\((1+hu^2w)^{m-1} 
	\le K(1+w)^{\max(m-1,0)}\),
\(\log(1+hu^2w)
	\le K+\log(1+w)\), and 
\(u^2w/(1+hu^2w)
	\le h_0^{-1}\).
Define \(\gamma=\min\{2m_0,m_0+1\}>1\). 
Since $2m-(m-1)_+\ge\gamma$, the density bounds imply
\begin{align*}
  |\psi_{m,h}(u,w)|
    +|\partial_m \psi_{m,h}(u,w)|
    +|\partial_h \psi_{m,h}(u,w)| 
  & \le K
      \begin{cases}
      w^{-1/2} & 0<w<1\\
      w^{-\gamma-1/2}(1+\log w) & 1\le w < \infty
      \end{cases}
\end{align*}
for some \(K < \infty\). The majorant is integrable and independent of $(m,h)$.

Differentiation under the integral yields the existence of both partial derivatives of \(\Psi\), while their joint continuity follows from dominated convergence. Since the compact rectangle was chosen arbitrarily, \(\Psi\) is continuously differentiable on \((1/2,\infty)\times(0,\infty)\).
\end{proof}

\section*{Acknowledgments}
The author gratefully acknowledges the support of the WISS 2025 project 'IDA-lab Salzburg' (20204-WISS/225/197-2019 and 20102-F1901166-KZP). 
The author used AI during the preparation of this manuscript for language editing, literature exploration and cross-checking, computation, and assistance in reviewing the exposition and technical arguments. All mathematical results, proofs, interpretations, \texttt{R} codes, and cited references were independently checked and verified by the author, who takes full responsibility for the contents of the paper.


\end{document}